\documentclass[11pt]{amsart}
\usepackage{amscd}
\usepackage{amsmath}
\usepackage{amssymb}
\usepackage{latexsym}
\usepackage{xypic}

\input xy
\xyoption{all}

\newtheorem{thm}{Theorem}[]
\newtheorem{prop}[thm]{Proposition}
\newtheorem{lem}[thm]{Lemma}
\newtheorem{lem-def}[thm]{Lemma-Definition}
\newtheorem{cor}[thm]{Corollary}

\theoremstyle{definition}

\numberwithin{equation}{section}

\newcommand{\Ga}{\Gamma}

\newcommand{\nc}{\newcommand}
\nc{\on}{\operatorname} \nc{\ct}{\check{\mathfrak t}}
\nc{\Z}{{\mathbb Z}}
\nc{\C}{{\mathbb C}} \nc{\pone}{{\mathbb P}^1} \nc{\pa}{\partial}
\nc{\F}{{\mathcal F}} \nc{\arr}{\rightarrow}
\nc{\larr}{\longrightarrow} \nc{\al}{\alpha} \nc{\ri}{\rangle}
\nc{\lef}{\langle} \nc{\W}{{\mathcal W}} \nc{\la}{\lambda}
\nc{\ep}{\epsilon} \nc{\su}{\widehat{{\mathfrak s}{\mathfrak
l}}_2} \nc{\sw}{{\mathfrak s}{\mathfrak l}} \nc{\g}{{\mathfrak g}}
\nc{\h}{{\mathfrak h}} \nc{\n}{{\mathfrak n}}
\nc{\N}{\widehat{\n}} \nc{\G}{\widehat{\g}} \nc{\De}{\Delta}
\nc{\gt}{\widetilde{\g}} \nc{\one}{{\mathbf 1}}
\nc{\z}{{\mathfrak Z}} \nc{\La}{\Lambda} \nc{\wt}{\widetilde}
\nc{\wh}{\widehat} \nc{\cri}{_{\kappa_c}} \nc{\kk}{k}
\nc{\sun}{\widehat{\sw}_N} \nc{\si}{\sigma} \nc{\el}{\ell}
\nc{\bi}{\bibitem} \nc{\om}{\omega} \nc{\ol}{\overline}
\nc{\ds}{\displaystyle} \nc{\dzz}{\frac{dz}{z}}
\nc{\Res}{\on{Res}} \nc{\mc}{\mathcal} \nc{\Cal}{\mathcal}
\nc{\bb}{{\mathfrak b}} \nc{\ot}{\otimes} \nc{\R}{{\mc R}}
\nc{\yy}{{\mc Y}} \nc{\ga}{\gamma}

\nc{\us}{\underset} \nc{\opl}{\oplus} \nc{\beq}{\begin{equation}}
\nc{\Fq}{{\mathcal F}} \nc{\Mq}{{\mathcal M}} \nc{\Rep}{\on{Rep}}
\nc{\sssec}{\subsubsection} \nc{\ssec}{\subsection}
\nc{\lan}{\langle} \nc{\ran}{\rangle}

\nc{\D}{\mathcal D} \nc{\Vect}{\on{Vect}} \nc{\ghat}{\G}
\nc{\T}{\mc T} \nc{\Tloc}{\T^\g_{\on{loc}}} \nc{\vac}{|0\ran}
\nc{\Wick}{{\mb :}} \nc{\mb}{\mathbf} \nc{\delz}{\partial_z}
\nc{\K}{{\cali K}} \nc{\cali}{\mathcal} \nc{\li}{\mathfrak l}
\nc{\lt}{\widetilde{\li}} \nc{\astar}{a^*} \nc{\cA}{{\mc A}}
\nc{\ka}{\kappa}

\nc{\OO}{{\mc O}} \nc{\AutO}{\on{Aut}\OO} \nc{\DerO}{\on{Der}\OO}
\nc{\DerpO}{\on{Der}_+\OO} \nc{\Au}{{\mc A}ut} \nc{\mf}{\mathfrak}
\nc{\V}{{\mathbb V}} \nc{\hh}{\wh{\h}}

\nc{\pp}{{\mathfrak p}} \nc{\mm}{{\mathfrak m}}
\nc{\rr}{{\mathfrak r}} \nc{\ket}{\rangle} \nc{\zz}{{\mathfrak z}}
\nc{\gr}{\on{gr}} \nc{\Spe}{\on{Spec}} \nc{\rv}{\crho}
\nc{\can}{\on{can}}
\nc{\CC}{{\mathcal C}} \nc{\Op}{\on{Op}_G(D)}
\nc{\MOp}{\on{MOp}_G(D)} \nc{\Db}{{\mathbb D}} \nc{\ww}{w}

\nc{\af}{{\mathbb A}^1} \nc{\bs}{\backslash} \nc{\laa}{(\la_i)}
\nc{\zn}{(z_i)}

\nc{\cla}{\check{\la}} \nc{\cmu}{\check{\mu}}
\nc{\crho}{\check{\rho}} \nc{\chal}{\check{\al}}
\nc{\cc}{{\mathfrak c}}

\nc{\M}{{\mathbb M}}

\nc{\ZZ}{{\mc Z}}

\nc{\UU}{{\mathbb U}}

\nc{\Conn}{\on{Conn}(\Omega^{\crho})}
\nc{\Con}{\on{Conn}(\Omega^{-\rho})}
\nc{\Co}{\on{Conn}(\Omega^{\rho})}

\nc{\ppart}{(\!(t)\!)}
\nc{\pparu}{(\!(u)\!)}
\nc{\ppartal}{(\!(t_\al)\!)}
\nc{\ppartinv}{(\!(t^{-1})\!)}
\nc{\zpart}{(\!(z)\!)}
\nc{\ppzi}{(\!(t-z_i)\!)} \nc{\ppinf}{(\!(t^{-1})\!)}
\nc{\Ind}{\on{Ind}} \nc{\I}{{\mathbb I}} \nc{\ppars}{(\!(s)\!)}
\nc{\QCoh}{\on{QCoh}}

\nc{\aff}{{\on{aff}}}
\nc{\AD}{{\mathbb A}}

\nc{\cG}{\check{G}}
\nc{\cB}{\check{B}}
\nc{\cT}{\check{T}}
\nc{\cg}{\check{\g}}
\nc{\cb}{\check{\bb}}
\nc{\ad}{\on{ad}}
\nc{\ab}{{\mathfrak a}}
\nc{\Loc}{\on{Loc}}
\nc{\Bun}{\on{Bun}}

\nc{\rk}{{\on{rk}(\cg)}}
\nc{\gm}{{\mathbb G}_m}
\nc{\cS}{\check{S}}
\nc{\ch}{\check{H}}
\nc{\Gm}{{\mathbb G}_m}

\begin{document}

\title{A rigid irregular connection on the projective line}

\author[Edward Frenkel]{Edward Frenkel$^1$}\thanks{$^1$Supported by
  DARPA and AFOSR through the grant FA9550-07-1-0543}

\address{Department of Mathematics, University of California,
Berkeley, CA 94720, USA}

\author{Benedict Gross}

\address{Department of Mathematics, Harvard University, Cambridge, MA
  02138, USA}

\dedicatory{To Victor Kac and Nick Katz on their 65th birthdays}

\date{January 2009; Revised April 2009}

\begin{abstract}

In this paper we construct a connection $\nabla$ on the trivial
$G$-bundle on $\pone$ for any simple complex algebraic group $G$,
which is regular outside of the points $0$ and $\infty$, has a regular
singularity at the point $0$, with principal unipotent monodromy, and
has an irregular singularity at the point $\infty$, with slope $1/h$,
the reciprocal of the Coxeter number of $G$.  The connection $\nabla$,
which admits the structure of an oper in the sense of Beilinson and
Drinfeld, appears to be the characteristic $0$ counterpart of a
hypothetical family of $\ell$-adic representations, which should
parametrize a specific automorphic representation under the global
Langlands correspondence. These $\ell$-adic representations, and their
characteristic $0$ counterparts, have been constructed in some cases
by Deligne and Katz. Our connection is constructed uniformly for any
simple algebraic group, and characterized using the formalism of
opers. It provides an example of the geometric Langlands
correspondence with wild ramification. We compute the de Rham
cohomology of our connection with values in a representation $V$ of
$G$, and describe the differential Galois group of $\nabla$ as a
subgroup of $G$.

\end{abstract}

\maketitle

\section{Introduction}

The Langlands correspondence relates automorphic representations of a
split reductive group $G$ over the ring of ad\`eles of a global field
$F$ and $\ell$-adic representations of the Galois group of $F$ with
values in a (slightly modified) dual group of $G$ (see Section
\ref{sect two}). On the other hand, the trace formula gives us an
effective tool to find the multiplicities of automorphic
representations satisfying certain local conditions.  In some cases
one finds that there is a unique irreducible automorphic
representation with prescribed local behavior at finitely many
places. A special case of this, analyzed in \cite{G3}, occurs when $F$
is the function field of the projective line $\pone$ over a finite
field $\kk$, and $G$ is a simple group over $\kk$.  We specify that
the local factor at one rational point of $\pone$ is the Steinberg
representation, the local factor at another rational point is a simple
supercuspidal representation constructed in \cite{GR}, and that the
local representations are unramified at all other places of $F$. In
this case the trace formula shows that there is an essentially unique
automorphic representation with these properties.  Hence the
corresponding family of $\ell$-adic representations of the Galois
group of $F$ to the dual group $\cG$ should also be unique.  An
interesting open problem is to find it.

Due to the compatibility of the local and global Langlands
conjectures, these $\ell$-adic representations should be unramified at
all points of $\pone$ except for two rational points $0$ and
$\infty$. At $0$ it should be tamely ramified, and the tame inertia
group should map to a subgroup of $\cG$ topologically generated by a
principal unipotent element. At $\infty$ it should be wildly ramified,
but in the mildest possible way.

Given a representation $V$ of the dual group $\cG$, we would obtain an
$\ell$-adic sheaf on $\pone$ (of rank $\dim V$) satisfying the same
properties. The desired lisse $\ell$-adic sheaves on ${\mathbb G}_m$
have been constructed by P. Deligne \cite{Trig} and N. Katz \cite{K}
in the cases when $\cG$ is $SL_n, Sp_{2n}, SO_{2n+1}$ or $G_2$ and $V$
is the irreducible representation of dimension $n, 2n, 2n+1$, and $7$,
respectively. However, there are no candidates for these $\ell$-adic
representations known for other groups $\cG$.

In order to gain a better understanding of the general case, we
consider an analogous problem in the framework of the geometric
Langlands correspondence. Here we switch from the function field $F$
of a curve defined over a finite field to an algebraic curve $X$ over
the complex field.  In the geometric correspondence (see, e.g.,
\cite{F:rev}) the role of an $\ell$-adic representation of the Galois
group of $F$ is played by a flat $\cG$-bundle on $X$ (that is, a pair
consisting of a principal $\cG$-bundle on $X$ and a connection
$\nabla$, which is automatically flat since $\dim X = 1$). Hence we
look for a flat $\cG$-bundle on $\pone$ having regular singularity at
a point $0 \in \pone$ with regular unipotent monodromy and an
irregular singularity at another point $\infty \in \pone$ with the
smallest possible slope $1/h$, where $h$ is the Coxeter number of
$\cG$ (see \cite{D} and Section \ref{oper connection} for the
definition of slope). By analogy with the characteristic $p$ case
discussed above, we expect that a flat bundle satisfying these
properties is unique (up to the action of the group ${\mathbb G}_m$ of
automorphisms of $\pone$ preserving the points $0,\infty$).

In this paper we construct this flat $\cG$-bundle for any simple
algebraic group $\cG$. A key point of our construction is that this
flat bundle is equipped with an oper structure. The notion of oper was
introduced by A. Beilinson and V. Drinfeld \cite{BD} (following the
earlier work \cite{DS}), and it plays an important role in the
geometric Langlands correspondence. An oper is a flat bundle with an
additional structure; namely, a reduction of the principal
$\cG$-bundle to a Borel subgroup $\check{B}$ which is in some sense
transverse to the connection $\nabla$. In our case, the principal
$\cG$-bundle on $\pone$ is actually trivial, and the oper
$\check{B}$-reduction is trivial as well. If $N$ is a principal
nilpotent element in the Lie algebra of a Borel subgroup opposite to
$\check{B}$ and $E$ is a basis vector of the highest root space for
$\check{B}$ on $\cg = \on{Lie}(\cG)$, then our connection takes the
form
\begin{equation}    \label{nabla1}
\nabla = d + N \frac{dt}{t} + E dt,
\end{equation}
where $t$ is a parameter on $\pone$ with a simple zero at $0$ and a
simple pole at $\infty$. We also give in Section \ref{twisted} a
twisted analogue of this formula, associated to an automorphism of
$\cG$ of finite order preserving $\cB$ (answering a question raised by
P. Deligne).

For any representation $V$ of $\cG$ our connection gives rise to a
flat connection on the trivial vector bundle of rank $\dim V$ on
$\pone$. We examine this connection more closely in the special cases
analyzed by Katz in \cite{K}. In these cases Katz constructs not only
the $\ell$-adic sheaves, but also their counterparts in characteristic
$0$, so we can compare with his results. These special cases share the
remarkable property that a regular unipotent element of $\cG$ has a
single Jordan block in the representation $V$. For this reason our
oper connection can be converted into a scalar differential operator
of order equal to $\dim V$ (this differential operator has the same
differential Galois group as the original connection). We compute this
operator in all of the above cases and find perfect agreement with the
differential operators constructed by Katz \cite{K}. This strongly
suggests that our connections are indeed the characteristic $0$
analogues of the special $\ell$-adic representations whose existence
is predicted by the Langlands correspondence and the trace formula.

\medskip

Another piece of evidence is the vanishing of the de Rham cohomology
of the intermediate extension to $\pone$ of the ${\mc D}$-module on
${\mathbb G}_m$ defined by our connection with values in the adjoint
representation of $\cG$. This matches the expectation that the
corresponding cohomology of the $\ell$-adic representations also
vanish, or equivalently, that the their global $L$-function with
respect to the adjoint representation of $\cG$ is equal to $1$. We
give two proofs of the vanishing of this de Rham cohomology. The first
uses non-trivial results about the principal Heisenberg subalgebras of
the affine Kac--Moody algebras due to V. Kac \cite{Kac1,Kac}. The
second uses an explicit description of the differential Galois group
of our connection and its inertia subgroups \cite{Katz:inv}.

Since the first de Rham cohomology is the space of infinitesimal
deformations of our local system (preserving its formal types at the
singular points $0$ and $\infty$) \cite{Katz2,BE,A}, its vanishing
means that our local system on $\pone$ is rigid.  We also prove the
vanishing of the de Rham cohomology for small representations
considered in \cite{K}. This is again in agreement with the vanishing
of the cohomology of the corresponding $\ell$-adic representations
shown by Katz. Using our description of the differential Galois group
of our connection and a formula of Deligne \cite{D} for the Euler
characteristic, we give a formula for the dimensions of the de Rham
cohomology groups for an arbitrary representation $V$ of $\cG$.

Finally, we describe some connections which are closely related to
$\nabla$, and others which are analogous to $\nabla$ coming from
subregular nilpotent elements. We also use $\nabla$ to give an example
of the geometric Langlands correspondence with wild ramification.

\medskip

The paper is organized as follows. In Section \ref{sect two} we
introduce the concepts and notation relevant to our discussion of
automorphic representations. In Section \ref{automorphic} we give the
formula for the multiplicity of automorphic representations from
\cite{G3}. This formula implies the existence of a particular
automorphic representation. In Section \ref{ell-adic} we summarize
what is known about the corresponding family of $\ell$-adic
representations. We then switch to characteristic $0$. In Section
\ref{oper connection} we give an explicit formula for our connection
for an arbitrary complex simple algebraic group, as well as its
twisted version. In Section \ref{special} we consider the special
cases of representations on which a regular unipotent element has a
single Jordan block. In these cases our connection can be represented
by a scalar differential operator. These operators agree with those
found earlier by Katz \cite{K}.

We then take up the question of computation of the de Rham cohomology
of our connection. After some preparatory material presented in
Sections \ref{coh}--\ref{nine} we prove vanishing of the de Rham
cohomology on the adjoint and small representations in Sections
\ref{van adj} and \ref{van small}, respectively. We also show that the
de Rham cohomology can be non-trivial for other representations using
the case of $SL_2$ as an example in Section \ref{sl2}. In Section
\ref{diff Galois} we determine the differential Galois group of our
connection. We then use it in Section \ref{dim of coh} to give a
formula for the dimensions of the de Rham cohomology for an arbitrary
finite-dimensional representation of $\cG$. In particular, we give an
alternative proof of the vanishing of de Rham cohomology for the
adjoint and small representations. In Section \ref{nearby} we discuss
some closely related connections.

Finally, in Section \ref{comments} we describe what the geometric
Langlands correspondence should look like for our connection.

\medskip

\noindent{\bf Acknowledgments.}  We met and started this project while
we were both visiting the Institut Math\'ematique de Jussieu in
Paris. E.F. thanks Fondation Sciences Math\'ematiques de Paris for its
support and the group ``Algebraic Analysis'' at Universit\'e Paris VI
for hospitality during his stay in Paris.

We have had several discussions of this problem with Dennis Gaitsgory
and Mark Reeder, and would like to thank them for their help. We thank
Dima Arinkin and Pierre Schapira for answering our questions about
${\mc D}$-modules and deformations, and Pierre Deligne for his help on
the setup of the global Langlands correspondence. We owe a particular
debt of gratitude to Nick Katz, who explained his beautiful results
carefully, and guided us in the right direction.

\section{Simple algebraic groups over global function fields}
\label{sect two}

Let $\kk$ be a finite field, of order $q$. Let $G$ be an absolutely
almost simple algebraic group over $\kk$ (which we will refer to as a
simple group for brevity). The group $G$ is quasi-split over $\kk$, and
we fix a maximal torus $A \subset B \subset G$ contained in a Borel
subgroup of $G$ over $\kk$. Let $\kk'$ be the splitting field of $G$,
which is the splitting field of the torus $A$, and put $\Gamma =
\on{Gal}(\kk'/\kk)$. Then $\Gamma$ is a finite cyclic group, of order
$1$, $2$, or $3$. Let $Z$ denote the center of $G$, which is a finite,
commutative group scheme over $\kk$.

Let $\cG$ denote the complex dual group of $G$. This comes with a
pinning $\check{T} \subset \check{B} \subset \cG$, as well as an
action of $\Gamma$ which permutes basis vectors $X_{-\al}$ of the
simple negative root spaces. The principal element $N = \sum X_{-\al}$
in $\check\g = \on{Lie}(\cG)$ is invariant under $\Gamma$
(\cite{G1}). Let $\check{Z}$ denote the finite center of $\cG$, which
also has an action of $\Gamma$.

There is an element $\ep$ in $Z(\cG)^\Gamma$ which satisfies $\ep^2=1$
and is defined as follows. Let $2\rho$ be the co-character of $\cT$
which is the sum of positive co-roots, and define
$$
\ep = (2\rho)(-1).
$$
Since the value of $\rho$ on any root is integral, $\ep$ lies in
$Z(\cG)$. It is also fixed by $\Gamma$. We have
$$
\ep = 1 \quad \longleftrightarrow \quad \rho \; \text{is a
  co-character of} \; \cT.
$$

In order to avoid choosing a square root of $q$ in the construction of
Galois representations, we will use the following modification
$\cG_1$ of $\cG$, which was suggested by Deligne. Let $\cG_1 = \cG
\times {\mathbb G}_m/(\ep \times -1)$. We have homomorphisms ${\mathbb
G}_m \to \cG_1 \to {\mathbb G}_m$ with composite $z \mapsto z^2$, The
group $\Ga$ acts on $\cG_1$, trivially on ${\mathbb G}_m$. The
co-character group $X_\bullet(\cT_1)$ contains the direct sum
$X_\bullet(\cT) \times X_\bullet({\mathbb G}_m)$ with index $2$. The
advantage of passing to $\cG_1$ is that we can choose a co-character
$$
\eta: {\mathbb G}_m \to \cT_1
$$
fixed by $\Ga$ which satisfies
$$
\langle \eta,\alpha \rangle = 1
$$
for all simple roots $\al$ of $\cG$. (This is impossible to do for
$\cG$ when $\ep \neq 1$.) Having chosen $\eta$, we let $w(\eta) \in
\Z$ be defined by composite map ${\mathbb G}_m \overset{\eta}\to \cG_1
\to {\mathbb G}_m$, $z \mapsto z^{w(\eta)}$. Then $w(\eta)$ is odd
precisely when $\ep \neq 1$.

\medskip

Let $X$ be a smooth, geometrically connected, complete algebraic curve
over $\kk$, of genus $g$. Let $F = \kk(X)$ be the global function field of
$X$. We fix two disjoint, non-empty sets $S, T$ of places $v$ of $F$,
and define the degrees
\begin{align*}
\on{deg}(S) &= \sum_{v \in S} \on{deg} v \; \geq \; 1, \\
\on{deg}(T) &= \sum_{v \in T} \on{deg} v \; \geq \; 1.
\end{align*}
A place $v$ of $F$ corresponds to a $\on{Gal}(\ol{\kk}/\kk)$-orbit on the
set of points $X(\ol{\kk})$. The degree $\on{deg} v$ of $v$ is the
cardinality of the orbit. Let
$$
M_G = \bigoplus_{d \geq 2} V_d(1-d)
$$
be the motive of the simple group $G$ over $F=\kk(X)$ (\cite{G2}). The
spaces $V_d$, of invariant polynomials of degree $d$, are all rational
representations of the finite, unramified quotient $\Gamma$ of
$\on{Gal}(F^s/F)$. The Artin $L$-function of $V=V_d$, relative to the
sets $S$ and $T$, is defined by
$$
L_{S,T}(V,s) = \prod_{v \not\in S} \on{det}(1-\on{Fr}_v q_v^{-s}|V)^{-1}
\prod_{v \in T} \on{det}(1-\on{Fr}_v q_v^{1-s}|V).
$$
Here $\on{Fr}_v=\on{Fr}^{\on{deg} v}$, where $\on{Fr}$ is the
Frobenius generator of $\Gamma$, $x \mapsto x^q$, and $q_v =
q^{\on{deg} v}$. This is known to be a polynomial of degree $\dim
V(2g-2+ \on{deg} S + \on{deg} T)$ in $q^{-s}$ with integral
coefficients and constant coefficient $1$ (\cite{W}). We define
$$
L_{S,T}(M_G) = \prod_{d \geq 2} L_{S,T}(V_d,1-d),
$$
which is a non-zero integer. In the next section, we will use the
integer $L_{S,T}(M_G)$ to study spaces of automorphic forms on $G$ over
$F$. We end this section with some examples.

Let $2=d_1,d_2,\ldots,d_{\on{rk}(G)}=h$ be the degrees of generators
of the algebra of invariant polynomials of the Weyl group, where
$\on{rk}(G)=\dim A$ is the rank of $G$ over the splitting field $\kk'$
and $h$ is the Coxeter number. If $G$ is split,
$$
L_{S,T}(M_G) = \prod_{i=1}^{\on{rk}(G)} \zeta_{S,T}(1-d_i),
$$
where $\zeta_{S,T}$ is the zeta-function of the curve $X-S$ relative
to $T$. Now assume that $G$ is not split, but that $G$ is not of type
$D_{2n}$. Then each $V_d$ has dimension $1$, $\Gamma$ has order $2$,
and $\Gamma$ acts non-trivially on $V_d$ if and only if $d$ is
odd. Hence
$$
L_{S,T}(M_G) = \prod_{d_i \; \on{even}} \zeta_{S,T}(1-d_i) \;
\prod_{d_i \; \on{odd}} L_{S,T}(\ep,1-d_i),
$$
where $\ep$ is the non-trivial quadratic character of $\Gamma$.

\section{Automorphic representations}    \label{automorphic}

Let ${\mathbb A}$ be the ring of ad\`eles of the function field
$F=\kk(X)$. Then $G(F)$ is a discrete subgroup, with finite co-volume, in
$G(\AD)$. Let $L$ denote the discrete spectrum, which is a
$G(\AD)$-submodule of $L^2(G(F) \bs G(\AD))$. Any irreducible
representation $\pi$ of $G(\AD)$ has finite multiplicity $m(\pi)$ in
$L$.

We will count the sum of multiplicities over irreducible
representations $\pi = \wh\otimes \pi_v$ with specified local
behavior. Specifically, for $v \not\in S \cup T$, we insist that
$\pi_v$ be an unramified irreducible representation of $G(F_v)$, in
the sense that the open compact subgroup $G(\OO_v)$ fixes a non-zero
vector in $\pi_v$. At places $v \in S$, we insist that $\pi_v$ is the
Steinberg representation of $G(F_v)$. Finally, at places $v \in T$, we
insist that $\pi_v$ is a simple supercuspidal representation of
$G(F_v)$, of the following type (cf. \cite{GR}). We let $\chi_v: P_v
\to \mu_p$ be a given affine generic character of a pro-$p$-Sylow
subgroup $P_v \subset G(\OO_v)$. We recall that $\chi_v$ is
non-trivial on the simple affine root spaces of the Frattini quotient
of $P_v$. (This is the affine analogue of a generic character of the
unipotent radical of a Borel subgroup.) Extend $\chi_v$ to a
character of $Z(q_v) \times P_v$ which is trivial on $Z(q_v)$; then
the compactly induced representation $\on{Ind}_{P_v \times
Z(q_v)}^{G(F_v)}(\chi_v)$ of $G(F_v)$ is multiplicity-free, with $\#
Z(\cG)^\Ga$ irreducible summands. We insist that $\pi_v$ be a
summand of this induced representation. 

The condition imposed on $\pi_v$ only depends on the
$I_v$-orbit of the generic character $\chi_v$, where $I_v$ is
the Iwahori subgroup of $G(F_v)$ which normalizes $P_v$. The
quotient group $I_v/P_v$ is isomorphic to $A(q_v)$, where $A$ is
the torus in the Borel subgroup, and there are
$(q_v-1) \cdot \# Z(q_v)$ different orbits. We note that affine generic
characters form a principal homogeneous space for the group
$A^{\on{ad}}(q_v) \times {\mathbb F}^\times_{q_v}$, where the latter
group maps a local parameter $t$ to $\lambda t$.

Using the simple trace formula, and assuming that some results of
Kottwitz \cite{Kt} on the vanishing of the local orbital integrals of
the Euler--Poincar\'e function on non-elliptic classes extend from
characteristic zero to characteristic $p$, we obtain the following
formula for multiplicities \cite{G3}:

\bigskip

\noindent {\em Assume that $p$ does not divide $\# Z(\cG)$. Then}
$$
\sum_{\pi \; {\on{as} \; \on{above}}} m(\pi) = L_{S,T}(M_G) \cdot
\frac{\# Z(\cG)^\Gamma}{\prod_{v \in S} \# Z(\cG)^{\Gamma_v}
    (-1)^{\on{rk}(G)_v}} \cdot \frac{\# Z(q)}{\prod_{v \in T} \# Z(q_v)
    (-1)^{\on{rk}(G)_v}},
$$
{\em where $\on{rk}(G)_v$ is the rank of $G$ over $F_v$.}

\bigskip

In the special case where $X=\pone$ has genus $0$, $S= \{ 0 \}$ and $T
= \{ \infty \}$, we have
$$
L_{S,T}(V_d,s) = 1, \qquad \on{for} \; \on{all} \; d,s.
$$
Hence
$$
L_{S,T}(M_G) = 1.
$$
Since $\Gamma_v=\Gamma$ and $q_v=q$ in this case, we find that
$$
\sum m(\pi) = 1.
$$
Hence there is a unique automorphic representation in the discrete
spectrum which is Steinberg at $0$, simple supercuspidal for a fixed
orbit of generic characters at $\infty$, and unramified
elsewhere. This global representation is defined over the field of
definition of $\pi_\infty$, which is a subfield of ${\mathbb
Q}(\mu_p)$. We will consider its Langlands parameter in the next
section.

\section{$\ell$-adic sheaves on ${\mathbb G}_m$}    \label{ell-adic}

Consider the automorphic representation $\pi$ of the simple group $G$
over $F=\kk(t)$, described at the end of Section \ref{automorphic}. This
representation is defined over the field ${\mathbb Q}(\mu_p)$. Its
local components $\pi_v$ are unramified irreducible representations of
$G(F_v)$, for all $v \neq 0,\infty$. The representation $\pi_0$ is the
Steinberg representation, and $\pi_\infty$ is a simple supercuspidal
representation.

Associated to $\pi$ and the choice of a co-character $\eta: {\mathbb
G}_m \to \cT_1 \subset \cG_1$, as described in Section \ref{sect two},
we should (conjecturally) have a global Langlands parameter. This will
be a continuous homomorphism
$$
\varphi_\la: \on{Gal}(F^s/F) \to \cG_1({\mathbb
  Q}(\mu_p)_\la)\rtimes \Gamma
$$
for every finite place $\la$ of ${\mathbb Q}_\ell(\mu_p)$ dividing a
rational prime $\ell \neq p$. Here we view $\cG_1$ as a pinned, split
group over ${\mathbb Q}$. The homomorphism $\varphi_\la$ should be
unramified outside of $0$ and $\infty$, and map $\on{Fr}_v$ to the
$\eta$-normalized Satake parameter for $\pi_v$ \cite{G:AMF}, which is
a semi-simple conjugacy class in $\cG_1({\mathbb Q}(\mu_p))\rtimes
\Gamma_v$. If this is true, the projection of $\varphi_\la$ to
${\mathbb G}_m({\mathbb Q}(\mu_p)_\la)$ will be the $w(\eta)$-power of
the $\ell$th cyclotomic character.

At $0$, $\varphi_\ell$ should be tamely ramified. The tame inertia
group should map to a $\Z_\ell$-subgroup of $\cG({\mathbb
Q}_\ell(\mu_p))$, topologically generated by a principal unipotent
element. At $\infty$, $\varphi_\ell$ should be wildly ramified, but
trivial on the subgroup $I_\infty^{1/h+}$ in the upper numbering
filtration, where $h$ is the Coxeter number of $G$.  If $p$ does not
divide $h$, the image of inertia should lie in the normalizer
$N(\check{T})$ of a maximal torus, and should have the form $E \cdot
\langle n \rangle$. Here $E \subset \cT[p]$ is a regular subgroup
isomorphic to the additive group of the finite field ${\mathbb
F}_p(\mu_h)$ and $n$ maps to a Coxeter element $w$ of order $h$ in the
quotient group $N(\check{T})/\check{T}$. The element $n$ is both
regular and semi-simple in $\cG$, and satisfies $n^h = \ep \in
Z(\cG)$, with $\ep^2 = 1$.

If
$$
\rho: \cG_1 \rtimes \Ga \to GL(V)
$$
is a representation over ${\mathbb Q}$, we would obtain from $\rho
\circ \varphi_\la$ a lisse $\la$-adic sheaf ${\mc F}$ on ${\mathbb
G}_m$ over $\kk$, with $\on{rank}({\mc F}) = \dim(V)$. Katz has
constructed and studied these lisse sheaves in the cases when the
Coxeter element in the Weyl group has a single orbit on the set of
non-zero weights for $\check{T}$ on $V$ (cf. Theorem 11.1 in
\cite{K}). In all of these cases, the principal nilpotent element $N =
\sum X_{-\al}$ has a single Jordan block on $V$. We list them in the
table below.

\bigskip

\begin{center}
\begin{tabular}{l|l}
$\cG$ & $\dim V$ \\
\hline
$SL_n$ & $n$ \\
\hline
$Sp_{2n}$ & $2n$ \\
\hline
$SO_{2n+1}$ & $2n+1$ \\
\hline
$G_2$ & $7$ \\
\end{tabular}
\end{center}

\bigskip

In all of these cases, there are no $I_\infty$-invariants on ${\mc F}$
and $sw_\infty({\mc F}) = 1$. If $j: {\mathbb G}_m \hookrightarrow
\pone$ is the inclusion, Katz has shown that $H^i(\pone,j_* {\mc F}) =
0$ for all $i$. Hence $L(\pi,V,s) = 1$.

More generally, consider the adjoint representation $\check\g$. Then
the sheaf ${\mc F}$ has rank equal to the dimension of $\cG$.  In this
case, we predict that

\begin{itemize}

\item[1)] $I_0$ has $\on{rk}(\cG)$ Jordan blocks on ${\mc F}$;

\item[2)] $I_\infty$ has no invariants on ${\mc F}$, and
  $sw_\infty({\mc F}) = \on{rk}(\cG)$;

\item[3)] $H^i(\pone,j_* {\mc F}) = 0$ for all $i$;

\item[4)] $L(\pi,\check\g,s)=1$.

\end{itemize}

Katz has verified this in the cases tabulated above, using the fact
that the adjoint representation $\check\g$  of $\cG$ appears in the
tensor product of the representation $V$ and its dual.

It is an open problem to construct the $\la$-adic sheaf ${\mc F}$ on
${\mathbb G}_m$ over the finite field $\kk$ for general groups $\cG$. We
consider an analogous problem, for local systems on ${\mathbb
G}_m(\C)= \C^\times$, in the next section.

\section{The connection}    \label{oper connection}

Let $\cG$ be a simple algebraic group over $\C$ and $\cg$ its Lie
algebra. Fix a Borel subgroup $\cB \subset \cG$ and a torus $\cT
\subset \cB$. For each simple root $\al_i$, we denote by $X_{-\al_i}$
a basis vector in the root subspace of $\cg = \on{Lie}(\cG)$
corresponding to $-\al_i$. Let $E$ be a basis vector in the root
subspace of $\cg$ corresponding to the maximal root $\theta$. Set
$$
N = \sum_{i=1}^\rk X_{-\al_i}.
$$

We define a connection $\nabla$ on the trivial $\cG$-bundle on $\pone$
by the formula
\begin{equation}    \label{nabla}
\nabla = d + N \frac{dt}{t} + E dt,
\end{equation}
where $t$ is a parameter on $\pone$ with a simple zero at $0$ and a
simple pole at $\infty$.

The connection $\nabla$ is clearly regular outside of the points $t=0$
and $\infty$, where the differential forms $\frac{dt}{t}$ and $dt$
have no poles.  We now discuss the behavior of \eqref{nabla} near the
points $t=0$ and $\infty$. It has regular singularity at $t=0$, with
the monodromy being a regular unipotent element of $\cG$. It also has
a double pole at $t=\infty$, so is irregular, and its slope there is
$1/h$, where $h$ is the Coxeter number. Here we adapt the definition
of slope from \cite{D}, Theorem 1.12: a connection on a principal
$\cG$-bundle with irregular singularity at a point $x$ on a curve $X$
has slope $a/b > 0$ at this point if the following holds. Let $s$ be a
uniformizing parameter at $x$, and pass to the extension given by
adjoining the $b^{\on{th}}$ root of $s$: $u^b = s$. Then the
connection, written using the parameter $u$ in the extension and a
particular trivialization of the bundle on the punctured disc at $x$
should have a pole of order $a + 1$ at $x$, and its polar part at $x$
should not be nilpotent.

To see that our connection has slope $1/h$ at the point $\infty$,
suppose first that $\cG$ has adjoint type. In terms of the
uniformizing parameter $s = t^{-1}$ at $\infty$, our connection has
the form
$$
d  - N \frac{ds}{s}  - E \frac{ds}{s^2}.
$$
Now take the covering given by $u^h = s$. Then the connection becomes
$$
d - h N \frac{du}{u} - h E \frac{du}{u^{h+1}}.
$$
If we now make a gauge transformation with $g = \rho(u)$ in the
torus $\cT$, where $\rho$ is the co-character of $\cT$ which is given by
half the sum of the positive co-roots for $\cB$, this becomes
\begin{equation}    \label{after}
d - h (N+E) \frac{du}{u^2} - \rho \frac{du}{u}.
\end{equation}
The element $N+E$ is regular and semi-simple, by Kostant's theorem
\cite{Ks}. Since the pole has order $(a+1) = 2$ with $a = 1$, the
slope is indeed $1/h$. If $\cG$ is not of adjoint type, then $g =
\rho(u)$ might not be in $\cT$, but it will be after we pass to
the covering obtained by extracting a square root of $u$. The
resulting slope will be the same.

Note that $\exp(\rho/h)$ is a Coxeter element in $\cG$, which
normalizes the maximal torus centralizing the regular element $N +
E$. We therefore have a close analogy between the local behavior of
our connection $\nabla$ and the desired local parameters in the
$\la$-adic representation $\varphi_\la$ at both zero and
infinity.

\bigskip

The connection we have defined looks deceptively simple. We now
describe how we used the theory of opers to find it. We recall from
\cite{BD} that a (regular) oper on a curve $U$ is a $\cG$-bundle with
a connection $\nabla$ and a reduction to the Borel subgroup $\cB$ such
that, with respect to any local trivialization of this
$\cB$-reduction, the connection has the form
\begin{equation}    \label{oper form}
\nabla = d + \sum_{i=1}^\rk \psi_i X_{-\al_i} + {\mathbf v},
\end{equation}
where the $\psi_i$ are nowhere vanishing one-forms on $U$ and ${\mathbf v}$
is a regular one-form taking values in the Borel subalgebra $\check\bb =
\on{Lie}(\cB)$. Here $X_{-\al_i}$ are non-zero vectors in the root
subspaces corresponding to the negative simple roots $-\al_i$ with
respect to any maximal torus $\cT \subset \cB$. Since the group $\cB$
acts transitively on the set of such tori, this definition is
independent of the choice of $\cT$. 

Any $\cB$-bundle on the curve $U = {\mathbb G}_m$ may be
trivialized. Therefore the space of $\cG$-opers on ${\mathbb G}_m$ may
be described very concretely as the quotient of the space of
connections \eqref{oper form}, where
$$
\psi_i \in \C[t,t^{-1}]^\times dt = \C^\times \cdot t^{\Z} dt,
$$
and ${\mathbf v} \in \check\bb[t,t^{-1}]dt$, modulo the gauge action
of $\cB[t,t^{-1}]$.

We now impose the following conditions at the points $0$ and $\infty$.
First we insist that at $t=0$ the connection has a pole of order $1$,
with principal unipotent monodromy. The corresponding
$\cB[t,t^{-1}]$-gauge equivalence class contains a representative
\eqref{oper form} with $\psi_i \in \C^\times t^{-1} dt$ for all
$i=1,\ldots,\rk$ and ${\mathbf v} \in \check\bb[t]dt$. By making a
gauge transformation by a suitable element in $\cT$, we can make
$\psi_i = \frac{dt}{t}$ for all $i$.

Second, we insist that the one-form ${\mathbf v}$ has a pole of order
$2$ at $t = \infty$, with leading term in the highest root space
$\cg_\theta$ of $\cb$, which is the minimal non-zero orbit for $\cB$
on its Lie algebra.  Then ${\mathbf v}(t) = Edt$, with $E$ a non-zero
vector of $\cg_\theta$.  These basis vectors are permuted
simply-transitively by the group ${\mathbb G}_m$ of automorphisms of
$\pone$ preserving $0$ and $\infty$ (that is, rescalings of the
coordinate $t$). Our oper connection therefore takes the form
\begin{equation}    \label{again}
\nabla = d + N \frac{dt}{t} + E dt.
\end{equation}
This shows that the oper satisfying the above conditions is unique up
to the automorphisms of $\pone$ preserving $0$ and $\infty$.

\medskip

We note that opers of this form (and possible additional regular
singularities at other points of $\pone$) have been considered in
\cite{FF:kdv}, where they were used to parametrize the spectra of
quantum KdV Hamiltonians. Opers of the form \eqref{again}, where $E$
is a regular element of $\cg$ (rather than nilpotent element
generating the maximal root subspace $\cg_\theta$, as discussed here),
have been considered in \cite{FFT,FFR}.  Finally, irregular
connections (not necessarily in oper form) with double poles and
regular semi-simple leading terms have been studied in \cite{B}.

\subsection{Twisted version}    \label{twisted}

In this paper we focus on the case of a ``constant'' group scheme over
${\mathbb G}_m$ corresponding to $\cG$. More generally, we may
consider group schemes over ${\mathbb G}_m$ twisted by automorphisms
of $\cG$. The connection $\nabla$ has analogues for these group
schemes which we now describe. It is natural to view them as complex
counterparts of the $\la$-adic representation $\varphi_\la$ for
non-split groups. This answers a question raised by Deligne.

Let $\sigma$ be an automorphism of $\cG$ of order $n$. It defines an
automorphism of the Lie algebra $\cg$, which we also denote by
$\sigma$. Define a group scheme $\cG_\sigma$ on $X=\Gm$ as
follows. Let $\wt{X} = \Gm$ be the $n$-sheeted cyclic cover of $X$
with a coordinate $z$ such that $z^n = t$. Then $\cG_\sigma$ is the
quotient of $\wt{X} \times \cG$, viewed as a group scheme over $X$, by
the automorphism $\wt\sigma$ which acts by the formula $(z,g) \mapsto
(ze^{2\pi i/n},\sigma(g))$. It is endowed with a (flat)
connection. Hence we have a natural notion of a $\cG_\sigma$-bundle on
$X$ with a (flat) connection. We now give an example of such an object
which is a twisted analogue of the flat bundles described above.

Take the trivial $\cG_\sigma$-bundle on $X$. Then a connection on it
may be described concretely as an operator
$$
\nabla = d + A(z) dz,
$$
where $A(z) dz$ is a $\wt\sigma$-invariant $\cg$-valued one-form on
$\wt{X}$.

Note that $\cG_\sigma$ and $\cG_{\sigma'}$ are isomorphic if $\sigma$
and $\sigma'$ differ by an inner automorphism of $\cG$. Hence, without
loss of generality, we may, and will, assume that $\sigma$ is an
automorphism of $\cG$ preserving the pinning we have chosen (thus,
$\sigma$ has order $1,2$ or $3$). In particular, it permutes the
elements $X_{-\al_i} \in \cg$. Then, for $A(z)$ as above, $z A(z)$
defines an element of the twisted affine Kac--Moody algebra
$\wh{\cg}_\sigma$, see \cite{Kac}. The lower nilpotent subalgebra of
this Lie algebra has generators $\wh{X}_{-\al_i},
i=0,\ldots,\ell_\sigma$, where $\ell_\sigma$ is the rank of the
$\sigma$-invariant Lie subalgebra of $\cg$, and $\wh{X}_{-\al_i},
i=1,\ldots,\ell_\sigma$, are $\sigma$-invariant linear combinations of
the elements $X_{-\al_j}$ of $\cg$ (viewed as elements of the constant
Lie subalgebra of $\wh{\cg}_\sigma)$. Explicit expressions for
$\wh{X}_{-\al_i}$ in terms of $\cg$ may be found in \cite{Kac}.

In the untwisted case the corresponding generators of the affine
Kac--Moody algebra (which is a central extension of $\cg[t,t^{-1}]$,
see Section \ref{van adj}) are $\wh{X}_{-\al_i}, i=i,\ldots,\ell =
\rk$, and $\wh{X}_{-\al_0} = Et$. Hence our connection \eqref{again}
may be written as
$$
\nabla = d + \sum_{i=0}^\ell \wh{X}_{-\al_i} \frac{dt}{t}.
$$

We now adapt the same formula in the twisted case and define the
following connection on the trivial $\cG_\sigma$-bundle on $\Gm$:
\begin{equation}    \label{twisted connection}
\nabla = d + \sum_{i=0}^{\ell_\sigma} \wh{X}_{-\al_i} \frac{dz}{z}.
\end{equation}
We propose \eqref{twisted connection} as a complex counterpart of the
$\la$-adic representation $\varphi_\la$ discussed above in the case
when the group $G$ is quasi-split and $\Gamma = \langle \sigma
\rangle$.

\medskip

The notion of oper may be generalized to the twisted case as
well. Namely, for any smooth complex curve $X$ equipped with an
unramified $n$-sheeted covering $\wt{X}$ we define the group scheme
$\cG_\sigma$ in the same way as above. Since $\sigma$ preserves a
Borel subgroup $\cB$ of $\cG$, the group scheme $\cG_\sigma$ contains
a group subscheme $\cB_\sigma$. The $\cG_\sigma$-opers on $X$
(relative to $\wt{X}$) are then $\cG_\sigma$-bundles on $X$ with a
connection $\nabla$ and a reduction to $\cB_\sigma$ satisfying the
condition that locally the connection has the form
\begin{equation}    \label{oper form1}
\nabla = d + \sum_{i=1}^{\ell_\sigma} \psi_i \wh{X}_{-\al_i} +
       {\mathbf v},
\end{equation}
where $\psi_i$ are nowhere vanishing one-forms and ${\mathbf v}$ takes
values in the Lie algebra of $\cB_\sigma$. (As above, if $\sigma$ acts
by permutation of elements $X_{-\al_j}$ of $\cg$, then the
$\wh{X}_{-\al_i}$ are $\sigma$-invariant linear combinations of the
$X_{-\al_j}$.)

It is clear that \eqref{twisted connection} is a $\cG_\sigma$-oper
connection on $\Gm$.

\section{Special cases}    \label{special}

If $V$ is any finite-dimensional complex representation of the group
$\cG$, a connection $\nabla$ on the principal $\cG$-bundle gives a
connection $\nabla(V)$ on the vector bundle ${\mc F}$ associated to
$V$. In our case, this connection is
\begin{equation}    
\nabla(V) = d + N(V) \frac{dt}{t} + E(V) dt,
\end{equation}
where $N(V)$ and $E(V)$ are the corresponding nilpotent endomorphisms
of $V$.

In this section, we will provide formulas for the first order matrix
differential operator $\nabla(V)_{t d/dt}$, for some simple
representations $V$ of $\cG$. We will be able to convert these matrix
differential operator into scalar differential operators, because in
these cases $N(V)$ will be represented by a principal nilpotent matrix
in $\on{End}(V)$. This will allow us to compare our connection with
the scalar differential operators studied by Katz in \cite{K}.

\medskip

{\bf Case I.} $\cG=SL_n$ and $V$ is the standard $n$-dimensional
representation with a basis of vectors $v_i, i=1,\ldots,n$, on which
the torus acts according to the weights $e_i$. Since $\al_i =
e_i-e_{i+1}, i=1,\ldots,n-1$ and $\theta=e_1-e_n$, we
can normalize the $X_{-\al_i}$ and $E$ in such a way that
\begin{equation}    \label{f e}
N(V) = \begin{pmatrix} 0 & 0 & 0 & ... & 0 & 0 \\
1 & 0 & 0 & ... & 0 & 0 \\
0 & 1 & 0 & ... & 0 & 0 \\
0 & 0 & 1 & ... & 0 & 0 \\
... & ... & ... & ... & ... & ... \\
0 & 0 & 0 & ... & 1 & 0
\end{pmatrix}, \qquad E(V) = \begin{pmatrix}
0 & 0 & 0 & ... & 0 & 1 \\ 
0 & 0 & 0 & ... & 0 & 0 \\
0 & 0 & 0 & ... & 0 & 0 \\
0 & 0 & 0 & ... & 0 & 0 \\
... & ... & ... & ... & ... & ... \\
0 & 0 & 0 & ... & 0 & 0
\end{pmatrix}.
\end{equation}

\smallskip

Therefore the operator $\nabla_{td/dt}$ of the connection
\eqref{nabla} has the form
\begin{equation}    \label{sl}
\nabla_{td/dt} = t\frac{d}{dt} +
\begin{pmatrix} 0 & 0 & 0 & ... & 0 & t \\
1 & 0 & 0 & ... & 0 & 0 \\
0 & 1 & 0 & ... & 0 & 0 \\
0 & 0 & 1 & ... & 0 & 0 \\
... & ... & ... & ... & ... & ... \\
0 & 0 & 0 & ... & 1 & 0
\end{pmatrix}
\end{equation}
and hence corresponds to the scalar differential operator
\begin{equation}    \label{scalar}
(t d/dt)^n + (-1)^{n+1} t.
\end{equation}

Now let $V$ be the dual of the standard representation. Then, choosing
as a basis the dual basis to the basis $\{ v_{n+1-i}
\}_{i=1,\ldots,n}$ (in the reverse order), we obtain the same matrix
differential operator \eqref{sl}. Hence the flat vector bundles
associated to the standard representation of $SL_n$ and its dual are
isomorphic.

\medskip

{\bf Case II.} $\cG=Sp_{2m}$ and $V$ is the standard $2m$-dimensional
representation. We choose the basis $v_i, i=1,\ldots,2m$, in which the
symplectic form is given by the formula
$$
\langle v_i,v_j \rangle = - \langle v_j,v_i \rangle =
\delta_{i,2m+1-j}, \qquad i<j.
$$

The weights of these vectors are
$e_1,e_2,\ldots,e_m,-e_m,\ldots,-e_2,-e_1$. Since
$\al_i = e_i-e_{i+1}, i=1,\ldots,m-1$, $\al_m=2e_m$ and
$\theta=2e_1$, we can normalize $N(V)$ and $E(V)$ in such a
way that they are given by formulas \eqref{f e}.

\smallskip

Hence $\nabla_{t d/dt}$ is also given by formula \eqref{sl} in this
case. It corresponds to the operator \eqref{scalar} with $n=2m$.

\medskip

{\bf Case III.} $\cG=SO_{2m+1}$ and $V$ is the standard
$(2m+1)$-dimensional representation of $SO_{2m+1}$ with the basis
$v_i, i=1,\ldots,2m+1$, in which the inner product has the form
$$
\langle v_i,v_j \rangle = (-1)^i \delta_{i,2m+2-j}.
$$
The weights of these vectors are
$e_1,e_2,\ldots,e_m,0,-e_m,\ldots,-e_2,-e_1$. Since
$\al_i = e_i-e_{i+1}, i=1,\ldots,m-1$, $\al_m=e_m$ and
$\theta=e_1+e_2$, we can normalize the $X_{-\al_i}$ and $E$ in
such a way that
\begin{equation}    \label{f e 1}
N(V) = \begin{pmatrix} 0 & 0 & 0 & ... & 0 & 0 \\
1 & 0 & 0 & ... & 0 & 0 \\
0 & 1 & 0 & ... & 0 & 0 \\
0 & 0 & 1 & ... & 0 & 0 \\
... & ... & ... & ... & ... & ... \\
0 & 0 & 0 & ... & 1 & 0
\end{pmatrix}, \qquad E(V) = \begin{pmatrix}
0 & 0 & 0 & ... & 1 & 0 \\ 
0 & 0 & 0 & ... & 0 & 1 \\
0 & 0 & 0 & ... & 0 & 0 \\
0 & 0 & 0 & ... & 0 & 0 \\
... & ... & ... & ... & ... & ... \\
0 & 0 & 0 & ... & 0 & 0
\end{pmatrix}.
\end{equation}

\smallskip

Therefore
\begin{equation}    \label{so}
\nabla_{td/dt} = t\frac{d}{dt} +
\begin{pmatrix} 0 & 0 & 0 & ... & t & 0 \\
1 & 0 & 0 & ... & 0 & t \\
0 & 1 & 0 & ... & 0 & 0 \\
0 & 0 & 1 & ... & 0 & 0 \\
... & ... & ... & ... & ... & ... \\
0 & 0 & 0 & ... & 1 & 0
\end{pmatrix}.
\end{equation}

We can convert this first order matrix differential operator into a
scalar differential operator. In order to do this, we need to find a
gauge transformation by an upper-triangular matrix which brings it to
a canonical form, in which we have $1$'s below the diagonal and other
non-zero entries occur only in the first row. This matrix is uniquely
determined by this property and is given by
$$
\begin{pmatrix} 1 & 0 & 0 & ... & t \\
0 & 1 & 0 & ... & 0 \\
0 & 0 & 1 & ... & 0 \\
... & ... & ... & ... & ... \\
0 & 0 & 0 & ... & 1
\end{pmatrix}.
$$
The resulting matrix operator is
$$
t\frac{d}{dt} +
\begin{pmatrix} 0 & 0 & 0 & ... & 2t & -t \\
1 & 0 & 0 & ... & 0 & 0 \\
0 & 1 & 0 & ... & 0 & 0 \\
0 & 0 & 1 & ... & 0 & 0 \\
... & ... & ... & ... & ... & ... \\
0 & 0 & 0 & ... & 1 & 0
\end{pmatrix}
$$
which corresponds to the scalar operator
$$
(td/dt)^{2m+1} - 2 t^2 d/dt - t.
$$

\medskip

{\bf Case IV.} $\cG=G_2$ and $V$ is the $7$-dimensional representation.

\smallskip

The Lie algebra $g_2$ of $G_2$ is a subalgebra of $so_7$. Both the
nilpotent elements $N$ and $E$ in $so_7$ may be simultaneously chosen
to lie in this subalgebra, where they are equal to the corresponding
elements for $g_2$.  Hence $\nabla_{td/dt}$ is equal to the operator
\eqref{so} with $m=3$, which corresponds to the scalar differential
operator
$$
(td/dt)^7 - 2 t^2 d/dt - t.
$$

\bigskip

There is one more case when $N$ is a regular nilpotent element in
$\on{End}(V)$; namely, when $\cG = SL_2$ and $V = \on{Sym}^n(\C^2)$ is
the irreducible representation of dimension $n+1$. We have already
considered the cases $n=1$ and $n=2$ (the latter corresponds to the
standard representation of $\cG=SO_3$).  These representations, and
the cases considered above, are the only cases when an oper may be
written as a scalar differential operator.

\medskip

The scalar differential operators we obtain agree with those
constructed by Katz in \cite{K}. More precisely, to obtain his
operators in the case of $SO_{2m+1}$ and $G_2$ we need to rescale $t$
by the formula $t \mapsto -\frac{1}{2} t$.

\medskip

In the above examples, the connection matrix was the same for
$Sp_{2n}$ as it was for $SL_{2n}$, and was the same for $G_2$ as it
was for $SO_7$.  The same phenomenon will occur for the pairs
$SO_{2n+1} < SO_{2n+2}, G_2 < D_4$, and $F_4 < E_6$, for the reason
explained in Section \ref{diff Galois}.

\section{De Rham cohomology}    \label{coh}

In the next sections of this paper, we will calculate the cohomology
of the intermediate extension of our local system to $\pone$, with
values in a representation $V$ of $\cG$. In particular, we will show
that this cohomology vanishes for the adjoint representation $\cg$, as
well as for the small representations tabulated in Section
\ref{special}. This is further evidence that our connection is the
characteristic $0$ analogue of the $\ell$-adic Langlands parameter. It
also implies that our connection is rigid, in the sense that it has no
infinitesimal deformations preserving the formal types at $0$ and
$\infty$, as such deformations form an affine space over the first de
Rham cohomology group \cite{Katz2,BE,A}.

\medskip

We begin with some general remarks on algebraic de Rham cohomology for
a principal $\cG$-bundle with connection $\nabla$ on the affine curve
$U ={\mathbb G}_m$. Any complex representation $V$ of $\cG$ then gives
rise to a flat vector bundle ${\mc F}(V)$ on $U$, where the connection
is $\nabla(V)$.

Since $U$ is affine, with the ring of functions $\C[t,t^{-1}]$, the
connection $\nabla(V)$ gives a $\C$-linear map
\begin{equation}    \label{C-lin}
\nabla(V): \quad V[t,t^{-1}] \to V[t,t^{-1}] \frac{dt}{t}.
\end{equation}
Any $\cG$-bundle on $U$ may be trivialized. Once we pick a
trivialization of our bundle ${\mc F}$, we represent the connection as
$\nabla = d + A$, where $A$ is a one-form on $U$ with values in the
Lie algebra $\cg$. Let $A(V)$ be the corresponding one-form on $U$
with values in $\on{End}(V)$. We may write
$$
A(V) = \sum_n A_n(V) \frac{dt}{t}
$$
with $A_n(V) \in \on{End}(V)$. If
$$
f(t) = \sum_n v_n t^n,
$$
we find that
$$
\nabla(V)(f) = \sum_n w_n t^n \frac{dt}{t}
$$
has coefficients $w_n$ given by the formula
$$
w_n = n v_n + \sum_{a+b=n} A_a(V)(v_b).
$$

For our specific connection
$$
\nabla = d + N(V) \frac{dt}{t} + E dt
$$
we find
$$
w_n = n v_n + N(V)(v_n) + E(V)(v_{n-1}).
$$

The ordinary de Rham cohomology groups $H^i(U,{\mc F}(V))$ are defined
as the cohomology of the complex \eqref{C-lin}:
\begin{align*}
H^0(U,{\mc F}(V)) &= \on{Ker} \nabla(V), \\
H^1(U,{\mc F}(V)) &= \on{Coker} \nabla(V).
\end{align*}
Thus elements of $H^0(U,{\mc F}(V))$ are solutions of the differential
equation $\nabla(V)(f) = 0$. For our particular connection, a solution
$f = \sum v_n t^n$ corresponds to a solution to the system of linear
equations
\begin{equation}    \label{sys lin}
n v_n + N(V)(v_n) + E(V)(v_{n-1}) = 0
\end{equation}
for all $n$. For example, if $v$ is in the kernel of both $N(V)$ and
$E(V)$, then $f=v$ is a constant solution, with $v_n=0$ for all $n
\neq 0$ and $v_0=v$.

We can also study the complex \eqref{C-lin} with functions and
one-forms on various subschemes of $U$. For example, the kernel and
cokernel of
\begin{equation}    \label{D0}
\nabla(V): \quad V\ppart \to V\ppart \frac{dt}{t}
\end{equation}
define the cohomology groups $H^0(D_0^\times,{\mc F}(V))$ and
$H^1(D_0^\times,{\mc F}(V))$, where $D^\times_0$ is the punctured disc
at $t=0$. Likewise, the kernel and cokernel of
\begin{equation}    \label{Dinfty}
\nabla(V): \quad V\ppartinv \to V\ppartinv \frac{dt}{t}
\end{equation}
define the cohomology groups $H^0(D_\infty^\times,{\mc F}(V))$ and
$H^1(D_\infty^\times,{\mc F}(V))$, where $D_\infty^\times$ is the
punctured disc at $t=\infty$. We will also identify the kernel and
cokernel of
\begin{equation}    \label{!coh}
\nabla(V): \quad V[[t,t^{-1}]] \to V[[t,t^{-1}]] \frac{dt}{t}
\end{equation}
as cohomology groups with compact support in Section \ref{nine}.

\medskip

The flat bundle ${\mc F}(V)$ on $U$ defines an algebraic, holonomic,
left ${\mc D}$-module on $U$ (which has the additional property of
being coherent as an ${\mc O}$-module). In the next section we will
recall the definition of the intermediate extension $j_{!*} {\mc
F}(V)$ in the category of algebraic, holonomic, left ${\mc
D}$-modules, where $j: U \hookrightarrow \pone$ is the inclusion. The
de Rham cohomology of $j_{!*} {\mc F}(V)$ may be calculated in terms
of some Ext groups in this category. We will establish the following
result, which is in agreement with the results of N. Katz on the
$\ell$-adic cohomology with coefficients in the adjoint representation
and small representations for the analogous $\ell$-adic
representations (in those cases in which they have been constructed).

\begin{thm}    \label{adjoint}
Assume that $\nabla = d + N \frac{dt}{t} + E dt$ and that $V$ is
either the adjoint representation $\cg$ of $\cG$ or one of the small
representations tabulated in Section \ref{special}. Then
$$
H^i(\pone,j_{!*} \mc F(V)) = 0
$$
for all $i$.
\end{thm}

We will provide two proofs of this result. The first, given in
Sections \ref{van adj} and \ref{van small}, uses the theory of affine
Kac--Moody algebras and the relation between the cohomology of the
intermediate extension of ${\mc F}(V)$ and solutions of the equation
$\nabla(V)(f) = 0$ in various spaces. The second, given in Section
\ref{dim of coh}, uses Deligne's Euler characteristic formula and a
calculation of the differential Galois group of $\nabla$. The latter
proof gives a formula for the dimensions of $H^i(\pone,j_{!*} {\mc
F}(V))$ for any representation $V$ of $\cG$.

\section{The intermediate extension and its cohomology}

Here we follow \cite{K}, Section 2.9 and \cite{BE} (see also
\cite{A}). Let $j: U \hookrightarrow \pone$ be the inclusion. We
consider the two functors
\begin{align}    \label{two fun}
j_* &= \text{direct image}, \\    \notag
j_! &= \Delta \circ j_* \circ \Delta
\end{align}
from the category of left holonomic ${\mc D}$-modules on $U$ to the
category of left holonomic ${\mc D}$-modules on $\pone$. The functor
$j_*$ is right adjoint to the inverse image functor, and $j_!$ is
defined using the duality functors $\Delta$ on these categories (see,
e.g. \cite{GM}, Section 5).

We have
$$
H^i(\pone,j_* {\mc F}) = H^i(U,{\mc F}), \qquad
H^i(\pone,j_! {\mc F}) = H^i_c(U,{\mc F}).
$$
The cohomology groups on $\pone$ are the Ext groups in the category of
holonomic ${\mc D}$-modules. The first equality follows from the
adjointness property of $j_*$, and the second can be taken as the
definition of the cohomology with compact support. Poincar\'e duality
gives a perfect pairing
$$
H^i(U,{\mc F}(V)) \times H^{2-i}_c(U,{\mc F}(V^*)) \to \C,
$$
where $V^*$ is the representation of $\cG$ that is dual to $V$. Thus,
we have $H^0_c(U,{\mc F}(V)) = 0$ and
\begin{equation}    \label{duality}
H^i_c(U,{\mc F}(V)) \simeq H^{2-i}(U,{\mc F}(V^*))^*, \qquad i=1,2.
\end{equation}

{}From the adjointness property of $j_*$, we obtain a map of ${\mc
  D}$-modules on $\pone$
$$
j_! {\mc F} \to j_* {\mc F}
$$
whose kernel and cokernels are ${\mc D}$-modules supported on $\{
0,\infty \}$. Let $j_{!*} {\mc F}$ be the image of $j_! {\mc F}$ in
$j_* {\mc F}$. We will now show that
\begin{align}    \label{zero}
H^0(\pone,j_{!*}{\mc F}(V)) &= H^0(U,{\mc F}(V)),
   \\ \label{one}
H^1(\pone,j_{!*}{\mc F}(V)) &= \on{Im}\left( H^1_c(U,{\mc F}(V)) \to
   H^1(U,{\mc F}(V)) \right), \\ \label{two}
H^2(\pone,j_{!*}{\mc F}(V)) &= H^2_c(U,{\mc F}(V)).
\end{align}
We will also describe an exact sequence involving the cohomology
groups on $D_0^\times$ and $D_\infty^\times$ which allows us to
compute $H^1(\pone,j_{!*}{\mc F}(V))$.

Let $t_\al$ be a uniformizing parameter at $\al = 0,\infty$ ($t$ and
$t^{-1}$, respectively) and let
$$
\delta_\al = \C\ppartal/C[[t_\al]]
$$
be the left delta ${\mc D}$-module supported at $\al$. We then have an
exact sequence of ${\mc D}$-modules on $\pone$
\begin{equation}    \label{first exact}
0 \to \bigoplus_\al H^0(D_\al^\times,{\mc F}) \otimes \delta_\al \to
j_! {\mc F} \to j_* {\mc F} \to \bigoplus_\al H^1(D_\al^\times,{\mc
  F}) \otimes \delta_\al \to 0.
\end{equation}

By the definition of $j_{!*} {\mc F}$, this gives two short exact
sequences
$$
0 \to \bigoplus_\al H^0(D_\al^\times,{\mc F}) \otimes \delta_\al \to
j_! {\mc F} \to j_{!*} {\mc F} \to 0
$$
$$
0 \to j_{!*} {\mc F} \to j_* {\mc F} \to \bigoplus_\al
H^1(D_\al^\times,{\mc   F}) \otimes \delta_\al \to 0
$$
We now take long exact sequence in cohomology and use the fact that 
\begin{align*}
H^0(\pone,\delta_\al) &= 0, \\
H^1(\pone,\delta_\al) &= \C.
\end{align*}
This gives a proof of \eqref{one}--\eqref{two}, and patching our two
long exact sequences along $H^1(\pone,j_{!*} {\mc F})$ gives a
six-term exact sequence
\begin{equation}    \label{six-term}
0 \to H^0(U,{\mc F}) \to \bigoplus_\al H^0(D_\al^\times,{\mc F}) \to
H^1_c(U,{\mc F}) \to
\end{equation}
$$
\to H^1(U,{\mc F}) \to \bigoplus_\al H^1(D_\al^\times,{\mc F}) \to
H^2_c(U,{\mc F}) \to 0.
$$
We will compare it later with an exact sequence obtained from the
snake lemma.

{}From the exact sequence \eqref{six-term} we deduce the following
condition for the vanishing of cohomology.

\begin{prop}    \label{van cond}
For a flat vector bundle ${\mc F}$ on $U$, we have $H^i(\pone,j_{!*}
\F) = 0$ for all $i$ if and only if

{\em (1)} $H^0(U,\F) = H^0(U,\F^*) = 0$;

{\em (2)} $\dim H^0(D_0^\times) + \dim H^0(D_\infty^\times) = \dim
   H^1_c(U,\F)$.
\end{prop}

\section{The dual complex}    \label{nine}

We have seen that the de Rham cohomology of the flat vector bundle
$\F(V)$ on $U$ can be calculated from the de Rham complex
\eqref{C-lin}. Since compactly supported cohomology of $\F(V)$ is dual
to the (ordinary) de Rham cohomology of $\F(V^*)$, it can be
calculated from the complex dual to
\begin{equation}    \label{C-lin*}
\nabla(V^*): \quad V^*[t,t^{-1}] \to V^*[t,t^{-1}] \frac{dt}{t}.
\end{equation}
In this section we will identify this dual complex with the complex
\begin{equation}    \label{C-lin-compl}
-\nabla(V): \quad V[[t,t^{-1}]] \to V[[t,t^{-1}]] \frac{dt}{t}
\end{equation}
by using the residue pairing at $t=0$, which is described in detail
below.

Hence $H^1_c(U,\F(V))$ is identified with the kernel of
\eqref{C-lin-compl} and $H^2_c(U,\F(V))$ with its cokernel. Using this
identification, we will compare the six-term exact sequence
\eqref{six-term} to the one obtained from the snake lemma.

We can also rewrite Proposition \ref{van cond} in a form that depends
only on solutions to $\nabla(f) = 0$.

\begin{cor}    \label{cor}
The cohomology of the intermediate extension $j_{!*} \F$ on $\pone$
vanishes if and only if

{\em (1)} $\on{Ker} \nabla(V) = 0$ on $V[t,t^{-1}]$ and $\on{Ker}
   \nabla(V^*) = 0$ on $V^*[t,t^{-1}]$;

{\em (2)} Every solution $f(t)$ to $\nabla(V)(f) = 0$ in
   $V[[t,t^{-1}]]$ can be written (uniquely) as a sum $f_0(t) +
   f_\infty(t)$, with $f_0$ and $f_\infty$ in $\on{Ker} \nabla(V)$ on
   $V\ppart$ and $V\ppartinv$, respectively.
\end{cor}

We now turn to the identification of \eqref{C-lin-compl} with the dual
of \eqref{C-lin*}. Define a bilinear pairing on $f \in V[[t,t^{-1}]]$
and $\omega \in V^*[t,t^{-1}] \frac{dt}{t}$ by
\begin{equation}    \label{res}
\langle f,\omega \rangle = \on{Res}_{t=0} 
S(f \cdot \omega),
\end{equation}
where $f \cdot \omega$ is the product in $V \otimes V^*[[t,t^{-1}]]
\frac{dt}{t}$ and $S: V \otimes V^* \to \C$ is the natural
contraction, so $S(f \cdot \omega)$ is an element of $\C[[t,t^{-1}]]
\frac{dt}{t}$. Explicitly, if $f = \sum v_n t^n$ and $\omega = \sum
\omega_m t^m \frac{dt}{t}$, then
$$
\langle f,\omega \rangle = \sum_{n+m=0} S(v_n \otimes \omega_m).
$$
This pairing identifies the direct product vector space
$V[[t,t^{-1}]]$ with the dual of the direct sum vector space
$V^*[t,t^{-1}] \frac{dt}{t}$. A similar pairing identifies 
$V^*[[t,t^{-1}]]$ with the dual of the direct sum vector space
$V[t,t^{-1}] \frac{dt}{t}$.

To complete the proof that this pairing identifies the dual of
\eqref{C-lin*} with \eqref{C-lin-compl}, we must show that the adjoint
of $\nabla(V^*)$ is $-\nabla(V)$. Write
$$
\nabla = d + \sum_m A_m t^m \frac{dt}{t},
$$
with $A_m \in \cg$. Then, for $g = \sum w_n t^n$ and $f$ as before we
have
\begin{align*}
\nabla(V^*)(g) &= dg + \sum_{m,n} A_m(V^*)(w_n) t^{m+n} \frac{dt}{t},
\\
\nabla(V)(f) &= df + \sum_{m,n} A_m(V)(v_n) t^{m+n} \frac{dt}{t}.
\end{align*}
The desired adjoint identity
$$
\langle \nabla(V)(f),g \rangle + \langle f,\nabla(V^*)(g) \rangle = 0
$$
then follows from the two identities
$$
\on{Res}_{t=0} (g \otimes df + f \otimes dg) = 0,
$$
$$
S(A(V)v,w) + S(v,A(V^*)w) = 0.
$$

\bigskip

We end this section with a reconstruction of the six-term exact
sequence \eqref{six-term}. The maps $\al(f) = (f,f)$ and $\beta(f,g) =
f-g$ give an exact sequence of vector spaces
$$
0 \to V[t,t^{-1}] \overset{\al}\to V\ppart \oplus V\ppartinv
\overset{\beta}\to V[[t,t^{-1}]] \to 0.
$$
Using $\nabla(V)$, we obtain a commutative diagram with exact rows
\begin{equation}    \label{diagram}
\begin{CD}
0 @>>> V[t,t^{-1}] @>{\al}>> V\ppart \oplus V\ppartinv
@>{\beta}>> V[[t,t^{-1}]] @>>> 0 \\ & &
@V{\nabla}VV  @V{(\nabla,\nabla)}VV  @V{\nabla}VV \\
0 @>>> V[t,t^{-1}] \frac{dt}{t} @>{\al}>> V\ppart \frac{dt}{t} \oplus
V\ppartinv \frac{dt}{t} @>{\beta}>> V[[t,t^{-1}]] \frac{dt}{t} @>>> 0
\end{CD}
\end{equation}
The $3$ kernels and the 3 cokernels have been identified with the
cohomology groups in \eqref{six-term}. Do the morphisms in
\eqref{six-term} come from an application of the snake lemma to
\eqref{diagram}? We note that we have made two sign choices: in the
pairing \eqref{res} we took the residue at $t=0$, not at $t=\infty$
(which would have changed the sign), and in the map $\beta$ we took
$f_0-f_\infty$, not $f_\infty-f_0$. We expect that with consistent
choice of signs the morphisms in \eqref{six-term} do indeed come from
\eqref{diagram} via the snake lemma.

\section{The vanishing of adjoint cohomology}    \label{van adj}

We now turn to the proof of Theorem \ref{adjoint}, using Corollary
\ref{cor}. Specifically, when $V$ is the adjoint representation of
$\cG$ or a small representation we will show that any solution $f(t) =
\sum v_n t^n$ of $\nabla(V)(f) = 0$ in $V[[t,t^{-1}]]$ satisfies
\begin{equation}    \label{neg}
v_n = 0 \qquad \text{for all} \quad n < 0.
\end{equation}
Next, we will use the following lemma.

\begin{lem}    \label{pos-neg}
Suppose that any solution $f = \sum v_n t^n \in V[[t,t^{-1}]]$ to
$\nabla(V)(f)=0$ satisfies property \eqref{neg} and the same property
holds if we replace $V$ by $V^*$. Then $H^i(\pone,j_{!*} \F(V)) = 0$
for all $i$.
\end{lem}

\begin{proof}
The equation $\nabla(V)(f)=0$ implies that the components $v_n$ satisfy
\begin{equation}    \label{rec}
n v_n + N(V)(v_n) + E(V)(v_{n-1}) = 0.
\end{equation}
If \eqref{neg} is satisfied, then it follows that $v_0$ lies in the
kernel of $N(V)$ on $V$. This also shows that there is a unique
solution
$$
f = \sum_{n \geq 0} v_n t^n
$$
for any $v_0$ in the kernel of $N(V)$, because $N(V)$ is nilpotent so
that the operator $n \on{Id} + N(V)$ is invertible on $V$ for all $n
\neq 0$. Clearly then, if $v_0 \neq 0$, this solution has
\begin{equation}    \label{pos}
v_n \neq 0  \qquad \text{for all} \quad n \geq 0.
\end{equation}
Hence there cannot be a non-zero solution to $\nabla(V)(f)=0$ that has
finitely many non-zero components for positive powers of $t$. We
obtain that
$$
H^0(U,\F(V)) = H^0(D_\infty^\times,\F(V)) = 0
$$
and
$$
H^0(D_0^\times,\F(V)) \simeq H^1_c(U,\F(V)) \simeq \on{Ker} N(V).
$$
Together we the same properties for $V$ replaced by $V^*$, this
implies that all of the criteria of Corollary \ref{cor} for the
vanishing of cohomology are met.
\end{proof}

We now turn to the proof of the property \eqref{neg} in the case when
$V=\cg$. We will drop $V$ in our notation and write $\nabla$ for
$\nabla(V)$, etc. The vector space $\cg[t,t^{-1}]$ is a $\Z$-graded
Lie algebra, with the Lie bracket
$$
[xt^n,yt^m] = [x,y]t^{n+m}.
$$
There is a $\cG$-invariant inner product
$$
\kappa: \cg \otimes \cg \to \C
$$
given by the Killing form (which is a unique such inner product up to
scaling). We define an inner product on $\cg[t,t^{-1}]$ by
\begin{equation}    \label{inner}
\left\langle \sum x_n t^n,\sum y_m t^m \right\rangle = \sum_{n+m=0}
\kappa(x_n,y_m).
\end{equation}
The $\Z$-grading on $\cg[t,t^{-1}]$ is given by the differential
operator $t d/dt$. If we write a solution $f(t)$ to $\nabla(f) =
0$ in terms of its graded pieces: $f = \sum v_n t^n$, then \eqref{rec}
becomes
\begin{equation}    \label{rec1}
n v_n + [N,v_n] + [E,v_{n-1}] = 0
\end{equation}
for all $n$.

To find the solutions to \eqref{rec1}, it is convenient to switch to a
different $\Z$-grading of $\cg[t,t^{-1}]$ called the principal
grading. Let $\rho$ be again half the sum of positive coroots for the
Borel subgroup $\cB$. Then
$$
[\rho,N] = -N, \qquad [\rho,E] = (h-1) E.
$$
The operator
\begin{equation}    \label{prin gr}
d = h \; t\frac{d}{dt} - \on{ad} \rho
\end{equation}
has integer eigenvalues on $\cg[t,t^{-1}]$, and defines the principal
grading with respect to which the element
\begin{equation}    \label{p1}
p_1 = N + Et
\end{equation}
has degree $1$. If we write a solution $f = \sum y_n$ of $\nabla(f)=0$
in its components for the principal grading, then \eqref{rec1} gives
rise to the identities
\begin{equation}    \label{hom comp}
n y_{n} + [\rho,y_n] + h [p_1,y_{n-1}] = 0
\end{equation}
for all $n$.

Since the eigenvalues of $\on{ad} \rho$ on $\cg$ are the integers in
the interval $[1-h,h-1]$, we see that the eigenvalues of $d$ on $\cg
  t^n$ are the integers of the form $nh + e$, with $1-h \leq e \leq
  h-1$. In particular, as eigenvector $y_m$ with eigenvalue $m$ has
  the form
$$
y_m \in \ct \; t^n,
$$
if $m=nh$, where $\ct \subset \cg$ is the unique Cartan subalgebra
containing $\rho$ (so $\ct$ is the kernel of $\on{ad} \rho$ on
$\cg$). If $n = nh + e$ with $0 < e < h$, then $y_m$ has the form
$$
y_m = a t^n + b t^{n+1}
$$
with $a$ of degree $e$ and $b$ of degree $(e-h)$  for $-\on{ad}
\rho$. From this one deduces that the component of $\cg[t,t^{-1}]$ of
degree $m$ with respect to the principal grading  has dimension equal
to the rank of $\cG$, except when $m$ is congruent to an exponent $e$
of $\cG$ (mod $h$), when the dimension is the rank of $\cG$ plus the
multiplicity of that exponent. Note that the original grading operator
$td/dt$ preserves the components with respect to the principal
grading.

V. Kac has studied the decomposition of $\cg[t,t^{-1}]$ under the
action of $\on{ad} p_1$, where $p_1=N+Et$. His results are summarized
in the following proposition. Set
$$
{\mathfrak a} = \on{Ker} (\on{ad} p_1), \qquad {\mathfrak c} = \on{Im}
    (\on{ad} p_{1}).
$$

\begin{prop}[\cite{Kac1}, Prop. 3.8]    \label{kac}
\hfill

{\em (1)} The Lie algebra $\cg[t,t^{-1}]$ has an orthogonal
    decomposition with respect to the inner product \eqref{inner},
    $$\cg[t,t^{-1}] = {\mathfrak a} \oplus {\mathfrak c}.$$

{\em (2)} ${\mathfrak a}$ is a commutative Lie subalgebra of
    $\cg[t,t^{-1}]$. With respect to the principal grading, $\ab =
    \bigoplus_{i \in I} \ab_i$, where $I$ is the set of all integers
    equal to the exponents of $\cg$ modulo the Coxeter number $h$, and
    $\dim \ab_i$ is equal to the multiplicity of the exponent $i \;
    \on{mod} \; h$.

{\em (3)} With respect to the principal grading, ${\mathfrak c} =
\bigoplus_{j \in \Z} {\mathfrak c}_j,$ where $\dim {\mathfrak c}_j =
\rk$, and the map $\on{ad} p_{1}: {\mathfrak c}_j \arr {\mathfrak
  c}_{j+1}$ is an isomorphism for all $j \in \Z$.
\end{prop}

Let $f$ be a solution to $\nabla(f)=0$ and $f=\sum y_n$ its
decomposition with respect to the principal grading. Then the
components satisfy \eqref{hom comp}. We now have the following crucial
lemma.

\begin{lem}    \label{zero1}
Suppose that the $y_n$ satisfy the equations \eqref{hom comp} and $y_n
\in \ab_n$ for some $n$. Then $y_m=0$ for all $m \leq n$.
\end{lem}

\begin{proof}
Applying \eqref{hom comp}, we obtain
\begin{equation}    \label{in cm}
t\frac{d}{dt} y_n = \frac{m}{h} y_n + \frac{1}{h}[\rho,y_n] \in
{\mathfrak c}_n,
\end{equation}
where ${\mathfrak c}_n$ is the degree $n$ homogeneous component of
${\mathfrak c} = \on{Im} (\on{ad} p_{1})$. Let us show that this is
impossible to satisfy if $y_n \in \ab_n$ and $y_n \neq 0$.

Consider the affine Kac--Moody algebra $\wh{\cg}$, which is the
universal central extension $\wh{\cg}$ of $\cg[t,t^{-1}]$ by
one-dimensional center spanned by an element ${\mathbf 1}$,
$$
0 \to \C{\mathbf 1} \to \wh{\cg} \to \cg[t,t^{-1}] \to 0
$$
The commutation relations in $\wh{\cg}$ read
$$
[At^n,Bt^m] = [A,B]t^{n+m} + n\kappa(A,B) \delta_{n,-m} {\mathbf 1}.
$$
According to \cite{Kac}, Lemma 14.4, the inverse image of $\ab$ in
$\wh{\cg}$ is a (non-degenerate) Heisenberg Lie subalgebra $\ab \oplus
\C {\mathbf 1}$. Hence there exists $z \in \ab_{-n}$ such that
$$
[y_n,z] \neq 0 \quad \on{in} \; \C {\mathbf 1} \subset \wh{\cg}.
$$
Write, for $n=kh+e$ (where $e$ is an exponent of $\cg$),
$$
y_n = a t^k + b t^{k+1}, \qquad z = a' t^{-k} +
b' t^{-k-1},
$$
as above, where $a, b, b, b'$ are homogeneous elements of
$\cg$ with respect to the grading defined by $-\ad \rho$ of degrees
$e, e-h, -e, h-e$, respectively. Then we find that
$$
[y_n,z] = (k \; \kappa(a,a') + (k+1) \; \kappa(b,b')) {\mathbf 1} \neq
0.
$$
But
$$
k \; \kappa(a,a') + (k+1) \; \kappa(b,b') = \langle t
\frac{d}{dt} \cdot y_n,z \rangle,
$$
where in the right hand side we use the inner product defined by
formula \eqref{inner}. This contradicts the condition that
$t\frac{d}{dt} \cdot y_n \in {\mathfrak c}_n$, which is orthogonal to
$\ab_{-n}$. Hence $y_n=0$.

Now, equation \eqref{hom comp} shows that if $y_n=0$, then
$y_{n-1} \in \ab_{n-1}$. Hence we find by induction that $y_m=0$ for
all $m \leq n$.
\end{proof}

Setting $n=0$ in \eqref{hom comp}, we obtain
$$
[\rho,y_0] + [p_1,y_{-1}] = 0.
$$
Since $\cg[t,t^{-1}]_0 = \ct = \on{Ker}(\on{ad} \rho)$, this shows
that
$$
[p_1,y_{-1}] = 0.
$$
Hence $y_{-1} \in \ab_{-1}$. Lemma \ref{zero2} then implies that
$y_n=0$ for all $n<0$. Thus, any solution $f \in \cg[[t,t^{-1}]]$ to
$\nabla(f)=0$ has the form $f=\sum_{n\geq 0} y_n$, and in particular
belongs to $\cg[[t]]$. Writing this solution as $f = \sum v_n t^n$, we
obtain that it satisfies property \eqref{neg}. Theorem \ref{adjoint}
for the adjoint representation now follows from Lemma \ref{pos-neg}
because $\cg^* = \cg$.

\section{Vanishing for small representations}    \label{van small}

Let $V$ be one of the small representations considered in Section
\ref{special}; that is, $n$-dimensional representation of $SL_n$,
$2n$-dimensional representation of $Sp_{2n}$, $(2n+1)$-dimensional
representation of $SO_{2n+1}$, or $7$-dimensional representation of
$G_2$. We denote the corresponding flat bundle on $U$ by ${\mc
F}(V)$. Since our connection for $Sp_{2n}$ and $G_2$ is the same as
for $SL_{2n}$ and $SO_7$, respectively, it is sufficient to consider
only the cases of $SL_n$ and $SO_{2n+1}$.

We will now prove Theorem \ref{adjoint} when $V$ is one of these
representations. We will follow the argument used in the proof of
Theorem \ref{adjoint} for the adjoint representation. Define a
$\Z$-grading on $V[t,t^{-1}]$ compatible with the principal
$\Z$-grading on $\cg[t,t^{-1}]$ defined by the operator $d$. The
representation $V$ has a basis $v_1,\ldots,v_p$ in which $N$ appears
as a lower Jordan block. We set
$$
\on{deg} v_i t^k = i-1 + kh.
$$
All graded components are one-dimensional in the case of $SL_n$. For
$SO_{2n+1}$ the components of degrees $kh, k \in \Z$, are
two-dimensional. Components of all other degrees are one-dimensional.
It is easy to see that the operator $td/dt$ preserves the graded
components.

We will use the operators $N(V)$ and $E(V)$ from Section
\ref{special}. Denote again $N(V) + E(V) t$ by $p_1$. Let $f \in
V[[t,t^{-1}]]$ be a solution to $\nabla(V)(f) = 0$. Decomposing $f =
\sum y_r$ with respect to the above grading, we obtain the following
system:
\begin{equation}    \label{hom comp1}
t \frac{d}{dt} y_{r} = - p_1 \cdot y_{r-1}, \qquad r \in \Z.
\end{equation}
The role of Lemma \ref{zero1} is now played by the following result.

\begin{lem}    \label{zero2}
Suppose that the $y_r$ satisfy the equations \eqref{hom comp1} and
$p_1 \cdot y_r = 0$ for some $r$. Then $y_m = 0$ for all $m\leq r$.
\end{lem}

\begin{proof}
In the case of $SL_n$ we have $\on{Ker} p_1 = 0$, so if $p_1 \cdot y_r
= 0$, then $y_r=0$. For $SO_{2n+1}$, $\on{Ker} p_1$ is spanned by the
vectors $v_1 t^k - v_{2n+1} t^{k-1}, k \in \Z$. Hence $y_r$ is one of
these vectors. Equation \eqref{hom comp1} tells us that
$$
t \frac{d}{dt}(v_1 t^k - v_{2n+1} t^{k-1}) = k v_1 t^{k} - (k-1)
v_{2n+1} t^{k-1}
$$
is in the image of $p_1$. But the image of $p_1$ in this
graded component is spanned by the vector $v_1 t^{k} + v_{2n+1}
t^{k-1}$. Since $k \in \Z$, it is impossible to satisfy this
condition. Therefore $y_r=0$.

Now, if $y_r =0$, then $p_1 \cdot y_{r-1} = 0$, by \eqref{hom
  comp1}. Hence we obtain by induction that $y_m = 0$ for all $m \leq
  r$.
\end{proof}

Now observe that $t d/dt$ annihilates the component of $V[t,t^{-1}]$
of degree $1$ (in fact, all components of degrees
$1,\ldots,h-1$). Hence we find from \eqref{hom comp1} with $r=1$ that
$p_1 \cdot y_0 = 0$. Therefore it follows from Lemma \ref{zero1} that
$y_m=0$ for all $m\leq 0$. Hence any solution $f = \sum f_n t^n$ to
$\nabla(V)(f)=0$ is a formal power series in $t$, that is, $f_n=0$ for
all $n<0$. Theorem \ref{adjoint} for small representations now follows
from Lemma \ref{pos-neg} and the fact that $V^* = V$ in the case of
$SO_{2n+1}$, and in the case of $SL_n$ the flat bundles associated to
$V$ and $V^*$ are isomorphic as well (see Case I in Section
\ref{special}).

\section{The case of $SL_2$}    \label{sl2}

In the previous two sections we have shown that the de Rham cohomology
of our connection vanishes in the adjoint representation as well as
the small representations. This is because solutions to the equation
$\nabla(V)(f) = 0$ enjoy special properties which do not hold for a
general representation $V$. The key property of $\cg$ and the small
representations that we have used is the fact that the weights of the
torus of a principal $SL_2$ on $V$ are all small (that is, the
eigenvalues of $\on{ad} \rho$ on $V$ have absolute value less than
$h$).

But for most other representations this is not the case and the de
Rham cohomology does not vanish. As an example, we consider in this
section the faithful irreducible representations $V = \on{Sym}^{2k-1}$
of $SL_2$ of even dimension $2k = n \geq 2$. We will see that there
are $k = n/2$ solutions in the space of formal power series in $t$ and
$t^{-1}$ with coefficients in $V$, but that only the zero solution
lies in $V\ppartinv$, and only one line of solutions lies in
$V\ppart$. Hence $H^1(\pone,j_{!*} \F(V))$ has dimension $k - 1 =
(n/2) - 1$.

We will see how to compute the dimensions of the de Rham cohomology of
$j_{!*} \F(V)$ on $\pone$, for any representation $V$ of $\cG$, in
Section \ref{diff Galois}.

\medskip

In the case of $SL_2$ our connection looks as follows:
$$
\nabla = d + F \frac{dt}{t} + E dt,
$$
where $F = X_{-\al_1}$ and $E$ are the standard generators of
$\sw_2$. Let $J=\on{diag}[1,2,\ldots,n]$. Make a change of
variables $t=z^2$ and apply gauge transformation by $z^{J}$. We
then obtain
$$
\nabla = d - J \frac{dz}{z} + 2(E+F) dz.
$$
Let $V_{\on{even}}$ (resp., $V_{\on{odd}}$) be the subspace of $V$ on
which the eigenvalues of $J$ are even (resp., odd).  The de Rham
complex on $U = {\mathbb G}_m$,
$$
V[t,t^{-1}] \overset{\nabla}\longrightarrow V[t,t^{-1}]dt,
$$
is identified with the subcomplex
\begin{equation}    \label{even}
V_{\on{even}}[z^2,z^{-2}] \oplus V_{\on{odd}}[z^2,z^{-2}] z
\overset{\nabla}\longrightarrow V_{\on{even}}[z^2,z^{-2}] \frac{dz}{z}
\oplus V_{\on{odd}}[z^2,z^{-2}] dz
\end{equation}
of the de Rham complex
\begin{equation}    \label{double}
V[z,z^{-1}] \overset{\nabla}\longrightarrow V[z,z^{-1}] dz
\end{equation}
on the double cover of $U$. Introduce a $\Z$-grading on this complex
by setting $\on{deg} v z^k = k$.

Note that the operator $E+F$ is conjugate to $2\rho =
\on{diag}[n-1,n-3,\ldots,-n+1]$. From now on we assume that $n$ is
even. Then the operator $E+F$ is invertible. Suppose that $y(z) = \sum
y_m z^m$ is in the kernel of \eqref{double}. Then we obtain the
following system of equations:
\begin{equation}    \label{recu}
(m \on{Id} - J) \cdot y_m = - 2(E+F) \cdot y_{m-1}
\end{equation}
on its homogeneous components. Let $m$ be the largest integer such
that $y_m = 0$ but $y_{m-1} \neq 0$. Then we should have $(E+F) \cdot
y_{m-1} = 0$, which is impossible. This shows that
$H^0(\pone,j_{!*}({\mc F}(V_n))) = 0$. Using duality, we obtain that
$H^2(\pone,j_{!*}({\mc F}(V_n))) = 0$.

In order to compute $H^1(\pone,j_{!*}({\mc F}(V_n))) = 0$, we first
compute the kernel of
\begin{equation}    \label{double1}
V[[z,z^{-1}]] \overset{\nabla}\longrightarrow V[[z,z^{-1}]] dz.
\end{equation}
Note that the operator $m \on{Id} - J$ is invertible for all $m\neq
1,\ldots,n$. Let us choose any $y_n \in V_n$. Then we can find $y_m,
m>n$, by using equation \eqref{recu} and inverting $m \on{Id} - J$,
and we can find $y_m, m<n$, by using \eqref{recu} and inverting
$E+F$. Thus, the kernel of \eqref{double1} is isomorphic to
$V_n$. An element of this kernel belongs to
$$
V_{\on{even}}[[z^2,z^{-2}]] \oplus V_{\on{odd}}[[z^2,z^{-2}]] z
$$
if and only if $y_n \in V_{\on{even}}$. Thus, we obtain that
$H^1(U,{\mc F}(V_n))^*$ is isomorphic to $V_{\on{even}}$. Using the
exact sequence \eqref{six-term}, we obtain that
$H^1(\pone,j_{!*}({\mc F}(V_n)))^*$ is the quotient of the above space
of solutions of \eqref{recu} by the subspace of those solutions for
which $y_m = 0$ for $m \gg 0$ or $m \ll 0$.

Those are precisely the solutions for which there exists
$m=1,\ldots,n$ such that $y_{m-1}=0$, but $y_m \neq 0$. We claim that
there are no such solutions for $m \neq n$. Indeed, denote by $v_i,
i=1,\ldots,n$, an eigenvector of $J$ with eigenvalue $i$. If
$y_{m-1}=0$, but $y_m \neq 0$, then $y_m$ is a multiple of $v_m$. But
$-2(E+F)(v_m)$ contains $v_{m+1}$ with non-zero coefficient if
$m<n$. Hence $-2(E+F)(v_m)$ cannot be in the image of $(m+1) \on{Id} -
J$, and so \eqref{recu} cannot be satisfied. If, on the other hand,
$y_n=v_n$, then $y_m=0$ for all $m<n$.

Thus, we obtain that there is a unique solution (up to a scalar) for
which $y_m = 0$ for $m \gg 0$ or $m \ll 0$. Hence $\dim
H^1(\pone,j_{!*}({\mc F}(V_n))) = (n/2) - 1$, which is non-zero if
$n\geq 4$.

\section{The differential Galois group}    \label{diff Galois}

In this section, we determine the differential Galois group of our
rigid irregular connection $\nabla$ on ${\mathbb G}_m$, as well as its
inertia subgroups (up to conjugacy) at $t=0$ and $t=\infty$.

We first review some of the general theory, which is due to N. Katz
\cite{Katz:inv}. Fix the point $t=1$ on ${\mathbb G}_m$. Then the
fiber at this point gives a fiber functor from the category of flat
complex algebraic vector bundles $(\F,\nabla)$ on ${\mathbb G}_m$ to
the category of finite-dimensional vector spaces. The automorphism
group of this fiber functor is, by definition, the differential Galois
group of ${\mathbb G}_m$. This is a pro-algebraic group over $\C$,
which we denote by $DG(\gm)$; Katz calls this group
$\pi_1^{\on{diff}}(\gm,1)$.

Since the fiber functor preserves tensor products, the fundamental
theorem of \cite{DM}, Ch. 2, gives an equivalence between the category
of flat bundles on $\gm$ with the category of finite-dimensional
representations of the differential Galois group $DG(\gm)$. If
$$
\varphi: DG(\gm) \to \on{GL}(V)
$$
corresponds to the flat bundle $(\F,\nabla)$, then
$$
H^0(\gm,\F) = \on{Ker} \nabla = V^{DG(\gm)}.
$$

Under this equivalence, a principal $\cG$-bundle with connection
$\nabla$ on $\gm$ defines a homomorphism
$$
\varphi_\nabla: DG(\gm) \to \cG
$$
up to conjugacy. The image $\cG_\nabla$ is an algebraic subgroup of
$\cG$, which we call the differential Galois group of $\nabla$.

At the two points $t=0$ and $t=\infty$ on $\pone - \gm$, we have local
inertia groups $I_0$ and $I_\infty$ in $DG(\gm)$, well-defined up to
conjugacy. Each inertia group $I=I_\al$ is filtered by normal
subgroups
$$
P^{x+} \subset P^x \subset P \subset I
$$
for rational $x>0$ (called the slopes). The wild inertia subgroup $P$
is a pro-torus over $\C$. As a pro-algebraic group over $\C$, $I/P$ is
isomorphic to the product of ${\mathbb G}_a$ with a pro-group $A$ of
multiplicative type, with character group $\C/\Z$. The additive part
comes from local systems with regular singularity at $\al$ with
unipotent monodromy, and the rest from the one-dimensional local
systems $d - a \; dt_\al/t_\al$ with solutions $t^a_\al$ near $t_\al =
0$, with $a \in \C/\Z$. The tame monodromy is then $e^{2\pi i a} \in
\C^\times$.  The torsion in the character group of $A$ is ${\mathbb
Q}/\Z$, so the component group of $A$ is the dual group $\wh\Z(1) =
\underset{\longleftarrow}\lim \; \mu_n$, which is the profinite Galois
group of $\C\ppartal$.

The quotient group $I/P$ acts on the pro-torus $P$ by
conjugation. Its connected component centralizes $P$; only the
component group $\wh\Z(1)$ acts non-trivially. On the character
group $\C$ of the quotient torus
$P^x/P^{x+}$, the element
$a = 2 \pi i m$ acts
by multiplication by the root of unity $e^{-ax}$. Hence the component
group acts through its finite quotient
$\wh\Z(1)/n\wh\Z(1) = \mu_n$, where $nx \equiv 0$ (mod $\Z$).

\bigskip

Now let
$$
\nabla = d + N \frac{dt}{t} + E dt
$$
be the connection introduced in Section \ref{oper connection}. We will
determine the differential Galois group $\cG_\nabla$ of $\cG$ by
studying the images $\varphi_\nabla(I_0)$ and $\varphi(I_\infty)$ in
$\cG$, and we begin with a discussion of the two homomorphisms
\begin{align*}
\varphi_\nabla: &I_0 \to \cG, \\
\varphi_\nabla: &I_\infty \to \cG
\end{align*}
up to conjugacy.

Since $\nabla$ has a regular singular point at $t=0$, the restriction
of $\varphi_\nabla$ to $I_0$ is trivial on the wild inertia subgroup
$P_0$. The resulting homomorphism is given by the analytic monodromy,
which maps $a=2\pi i n$ in $Z(1)$ to $\exp(-aN)$ in $\cG$. In
particular, the image $\varphi_\nabla(I_0)$ is an additive subgroup
${\mathbb G}_a = \exp(zN)$ of $\cG$, containing the principal
unipotent element $u = \exp(-2\pi i N)$.

To determine the image of the inertia group at $t = \infty$, we first
make the assumption that $\rho$ is a co-character of $\cG$. Let $u^h =
s = t^{-1}$. Then by our earlier results in Section \ref{oper
connection}, over the extension $\C\pparu$ of $\C\ppars$ our
connection is equivalent to
$$
d - h (N+E) \frac{du}{u^2} - \rho \frac{du}{u}.
$$
By a fundamental result of Kostant \cite{Ks}, the element $(N + E)$ is
regular and semi-simple. Hence the highest order polar term of our
connection over $\C((u))$ is diagonalizable, in any representation $V$
of $\cG$. It follows that the slopes of this connection over
$\C\pparu$, as defined by Deligne \cite{D}, Theorem 1.12, are either
$0$ or $1$, the former occurring at the zero eigenspaces for $N+E$ on
$V$ and the latter occurring at the non-zero eigenspaces. Since Katz
has shown \cite{Katz:inv}, Sections 1 and 2.2.11.2, that the slopes
over the extension $\C\pparu$ are $h$ times the slopes over the
original completion $\C\ppars$, we see that the original slopes are
either $0$ or $1/h$. In particular, $I_\infty$ is trivial on the
subgroup $P^{1/h+}$. A similar argument works when $\rho$ is not a
co-character, using the extension of degree $2h$.

The image $\cS := \varphi(P_\infty)$ of wild inertia subgroup
$P_\infty$ is the smallest torus in $\cG$, whose Lie algebra contains
the regular, semi-simple element $N+E$. The irregularity
$\on{Irr}_\al(V)$ of a representation $V$ of $\cG$ at an irregular
singular point $\al$ was defined by Deligne \cite{D}, p. 110, and
shown by Katz \cite{Katz:inv}, Sections 1 and 2.3, to be the sum of the
slopes. From the above analysis, we deduce that
$$
h \on{Irr}_\infty(V) = \dim V - \dim V^{\cS}.
$$

The full image $\ch = \varphi(I_\infty)$ normalizes $\cS$, and
the quotient is generated by the element $n=(2\rho)(e^{\pi i/h})$.
The element $n$ is regular and semi-simple in $\cG$. Further, $\ep: =
n^h$ is a central involution in $\cG$ which is equal to identity if
and only if $\rho$ is a co-character of $\cG$. The element $n$ acts
on $N+E$ by multiplication by $e^{-2\pi i/h}$. The irreducibility of
the $h^{\on{th}}$ cyclotomic polynomial over $\Z$ shows that $\cS$ has
dimension $\phi(h) = \#(\Z/h\Z)^\times$ and that the eigenvalues of
$n$ on $\on{Lie}(\cS)$ are the primitive $h^{\on{th}}$ roots of unity.

Since $n$ normalizes $\cS$, it also normalizes the centralizer $\cT'$
of $\cS$ in $\cG$, which is a maximal torus (note that it is different
from the maximal torus $\cT$ considered in Section \ref{oper
  connection}). The image $w$ of $n$ in $N(\cT')/\cT'$ is a Coxeter
class in the Weyl group, and has order $h$.

The character group $X^*(\cS)$ is the free quotient of $X^*(\cT')$
where $w$ acts by primitive $h^{\on{th}}$ roots of unity. Thus the
characters of $\cT'$ which restrict to the trivial character of $\cS$
are generated by those $\la: \cT' \to \gm$ where the $\langle w
\rangle$-orbit of $\la$ has size less than $h$.

This completes the description of the local inertia groups, and we now
turn to the global differential Galois group $\cG_\nabla$.

\begin{prop}    \label{contains}
Let $\cG_0 \subset \cG$ be an algebraic subgroup which contains
$\varphi(I_\infty) = \ch$ and $\varphi(I_0) = \exp(zN)$. Then
$\cG_0$ is reductive and contains the image of a principal $SL_2$ in
$\cG$.
\end{prop}

\begin{proof} Let $R(\cg_0)$ be the unipotent radical of $\cG_0$, and
  let $Z = \on{Lie}(Z)$ be its center. We will show that $Z=0$.

The group $\ch = \varphi(I_\infty)$ acts on $R(\cG_0)$ and $Z$. Since
$\cS$ contains regular elements, every root $\al: \cT' \to \gm$
restricts to a non-trivial character of $\cS$. Hence the action of
$\ch = \langle \cS,n \rangle$ on $\cg$ decomposes as
$$
\on{Lie}(\cT') \oplus \bigoplus_{i=1}^{\rk} W_i,
$$
where the $W_i$ are irreducible representations of dimension $h$ whose
restriction to $\cS$ contains an entire $w$-orbit of roots.

Since $Z$ is nilpotent and $\on{Lie}(\cT')$ is semi-simple, $Z$ must
be the sum of certain $W_i$. But each $W_i$ contains a non-zero vector
$v_0 = \sum_{i=1}^h n^i(v)$ fixed by $\langle n \rangle$ (as $n^h =
\ep$ is central and acts trivially on $\cg$). Since $n$ is regular and
semi-simple, $v_0$ is a semi-simple element and hence cannot be
contained in $Z$. Therefore $Z=0$.

Since $\cG_0$ is reductive and contains the principal unipotent
element $u$, it contains a principal $SL_2$ \cite{Carter}. We note
that a principal embedding $SL_2 \to \cG$ maps the central element $-1
\in SL_2$ to the element $\ep \in \ch \subset \cG$.
\end{proof}

Proposition \ref{contains} is a serious constraint on the global image
$\cG_\nabla$, as the reductive subgroups of $\cG$ containing a
principal $SL_2$ are severely limited by a result of \cite{SS} which
goes back to the work of Dynkin \cite{Dynkin}. They are all simple,
and their Lie algebras appear in one of the following maximal chains:
$$
\xymatrix{
\sw_2 \ar[r] & {\mathfrak s}{\mathfrak p}_{2n} \ar[r] & \sw_{2n} &
}
$$
$$
\xymatrix{
& & \sw_{2n+1} & \\
\sw_2 \ar[r] & {\mathfrak s}{\mathfrak o}_{2n+1} \ar[ur] \ar[dr] \\
& & {\mathfrak s}{\mathfrak o}_{2n+2} & \\
}
$$
$$
\xymatrix{& & & \sw_7 & \\
\sw_2 \ar[r] & {\mathfrak g}_2 \ar[r] & {\mathfrak s}{\mathfrak o}_7
\ar[ur] \ar[dr] \\
& & & {\mathfrak s}{\mathfrak o}_8 \\
\sw_2 \ar[r] & {\mathfrak f}_4 \ar[r] & {\mathfrak e}_6 & \\
\sw_2 \ar[r] & {\mathfrak e}_7 & & \\
\sw_2 \ar[r] & {\mathfrak e}_8 & & \\
}
$$

In some of these cases, the subgroup $\cG_0$ cannot contain $\ch$, as
its Coxeter number is less than the Coxeter number of $\cG$. Looking
at the minimal cases remaining, we obtain

\begin{cor}    \label{descr}
If $\cG$ is simple of type $A_{2n}, n \geq 1, C_n, n \geq 1, B_n, n
\geq 4$, $G_2$, $F_4$, $E_7$, $E_8$, then $\cG_\nabla = \cG$.
\end{cor}

Indeed, by the above list of embeddings and Proposition
\ref{contains}, any $\cG_0 \subset \cG$ containing $\ch$ and $\langle
n \rangle$ must be equal to $\cG$.

For the remaining cases, observe that the automorphism group $\Sigma$
of the pinning of $\cG$ is known to be isomorphic to the outer
automorphism group of $\cG$. This finite group fixes $N$ and acts on
the highest root space. If $\cG$ is not of type $A_{2n}$, it also
fixes $E$.  Hence $\Sigma$ fixes the connection $\nabla$, and its
differential Galois group $\cG_\nabla$ is contained in
$\cG^\Sigma$. Thus Corollary \ref{descr} gives the differential Galois
group in all cases.

\begin{cor}
If $\cG$ is of type $A_{2n-1}$, then $\cG_\nabla$ is of type $C_n$,
with center the kernel of the center of $\cG$ on the second exterior
power representation.

If $\cG$ is of type $D_{2n+1}$ with $n \geq 4$, then $\cG_\nabla$ is
of type $B_n$, with the center the kernel of the center of $\cG$ on
the standard representation.

If $\cG$ is of type $D_4$ or $B_3$, then $\cG_\nabla$ is of type
$G_2$.

If $\cG$ is of type $E_6$, then $\cG_\nabla$ is of type $F_4$.
\end{cor}

\section{The dimension of cohomology}    \label{dim of coh}

We now use the calculation of the differential Galois group
$\cG_\nabla$ of our connection, with its inertia subgroups at $0$ and
$\infty$, to determine the dimensions of the cohomology groups of the
${\mc D}$-modules $j_* \F(V)$, $j_! \F(V)$ and $j_{!*} \F(V)$ on
$\pone$ associated to a representation $V$ of $\cG$.

We will assume that we are in one of the cases described in Corollary
\ref{descr}. Otherwise, the computation of cohomology reduces to one
in a smaller group given by Corollary \ref{descr}. We will also assume
that $V$ is irreducible  non-trivial representation of $\cG$, so
\begin{align*}
H^0(\pone,j_* \F(V)) &= H^0(U,\F(V)) = V^{\cG} = 0, \\
H^2(\pone,j_! \F(V)) &= H^2_c(U,\F(V)) = \on{Hom}_{\cG}(V^*,\C) = 0.
\end{align*}
We will use Deligne's formula \cite{D}, Section 6.21.1, for the Euler
characteristic
$$
\chi(H^\bullet(U,\F)) = \chi(H^\bullet_c(U,\F)) = \chi(U)
\on{rank}(\F) - \sum_\al \on{Irr}_\al(\F).
$$
In our case, $\chi(U) = 0$ and $\nabla$ is regular at $\al=0$. Hence
$$
\on{dim} H^1(\pone,j_* \F) = \on{dim} H^1(\pone,j_! \F) =
\on{Irr}_\infty(\F).
$$

The kernel of the map $H^1_c(U,\F) \to H^1(U,\F)$ is isomorphic to the
direct sum
$$
H^0(D_0^\times,\F) \oplus H^0(D_\infty^\times,\F) = V^{I_0} \oplus
V^{I_\infty}
$$
by \eqref{six-term}. Hence we obtain

\begin{prop}
If $\cG_\nabla = \cG$ and $V$ is an irreducible, non-trivial
representation of $\cG$ with associated flat vector bundle $\F(V)$ on
$\gm$, then
$$
H^0(\pone,j_{!*} \F(V)) = H^2(\pone,j_{!*} \F(V)) = 0,
$$
\begin{equation}    \label{dV}
d(V) := \dim H^1(\pone,j_{!*} \F(V)) = \on{Irr}_\infty(V) - \dim V^{I_0}
- \dim V^{I_0}.
\end{equation}
\end{prop}

If we don't assume that $\cG_\nabla = \cG$ or that $V$ is an
irreducible, non-trivial representation of $\cG$, we obtain the
formulas
$$
\dim H^0(\pone,j_{!*} \F(V)) = \dim H^2(\pone,j_{!*} \F(V)) = \dim
V^{\cG_\nabla},
$$
$$
\dim H^1(\pone,j_{!*} \F(V)) =  \on{Irr}_\infty(V) - \dim V^{I_0}
- \dim V^{I_0} + 2 \dim V^{\cG_\nabla}.
$$

Since
$$
h  \on{Irr}_\infty (V) = \dim V - \dim V^{\cS},
$$
as we have seen above, this allows us to compute the dimensions of the
cohomology of the middle extension for any representation $V$ of
$\cG$, provided that we know the restriction of $V$ to the three
subgroups $\cS$, $\ch$ and $\cG_\nabla$, and the restriction of $V$
to a principal $SL_2$. We will now make this more explicit.

\bigskip

The irreducible representations of $SL_2$ all have the form
$\on{Sym}^k, k \geq 0$. Hence we may write the restriction of $V$ to
the principal $SL_2$ as
$$
V = \bigoplus_{k \geq 0} \; (\on{Sym}^k)^{\oplus m(k)}
$$
with multiplicities $m(k) \geq 0$. We then have
\begin{equation}    \label{I0}
\dim V^{I_0} = \sum_{k \geq 0} m(k)
\end{equation}
as $\varphi(I_0) = \exp(tN)$ fixes a unique line on each
irreducible factor. We note that the parity of these $k$ is
determined by $V$: if $m(k) > 0$ we have
$$
(-1)^k = \ep|_{V}.
$$

The irreducible complex representations of $\ch$ have dimensions
either $h$ or $1$. The irreducible $W$ of dimension $h$ restrict to a
sum of $h$ distinct non-trivial characters $\la$ of $\cS$, in a single
$\langle n \rangle$-orbit. The irreducible $\chi$ of dimension $1$ are
the representations trivial on $\cS$, and determined by $\chi(n)$,
which lies in $\mu_{2h}$. We may therefore write
\begin{equation}    \label{decomposition}
V = \bigoplus_i \chi_i^{\oplus m(\chi_i)} \oplus \bigoplus_j
W_j^{\oplus m(W_j)}.
\end{equation}
If $m_i>0$, then
$$
\chi_i(n)^h = \ep|_V.
$$
In terms of the decomposition \eqref{decomposition}, we have
\begin{equation}    \label{Iinfty}
\dim V^{I_\infty} = m(\chi_0),
\end{equation}
where $\chi_0$ is the trivial character, $\chi_0(n) = 1$. Each
one-dimensional representation of $\ch$ is tame, and each irreducible
$h$-dimensional representation has irregularity $1 = h \cdot
(1/h)$. Hence
\begin{align}    \label{Irr}
\on{Irr}_\infty(V) &= \sum_j m(W_j) \\ \notag
&= \frac{1}{h} \cdot \# \{ \text{non-trivial weight spaces for $\cS$
on $V$} \}.
\end{align}

If we know the restriction of $V$ to a principal $SL_2$ and to $\ch$,
then formulas \eqref{dV}, \eqref{I0}, \eqref{Iinfty}, and \eqref{Irr}
allow us to determine the dimension $d(V)$ of $H^1(\pone,j_{!*}
\F(V))$. For example, when $V=\cg$ is the adjoint representation, we
have $\ep=+1$ on $V$ and
\begin{align*}
V &= \bigoplus_{i=1}^{\rk} \on{Sym}^{2d_i-2}, \\
&= \bigoplus_{i=1}^{\rk} \chi_i \oplus \bigoplus_{j=1}^{\rk} W_j,
\end{align*}
where the $d_i$ are the degrees of invariant polynomials and
$\chi_i\neq 1$ for all $i$. Hence $\on{Irr}_\infty = \dim V^{I_0} =
\rk$, $V^{I_\infty} = 0$, and so $d(V)=0$. This gives the second proof
of the rigidity of our connection (Theorem \ref{adjoint}), for the
groups $\cG$ in Corollary \ref{descr}. (When $\cG_\nabla$ is a proper
subgroup of $\cG$, we find that $\cg = \cg_\nabla \oplus V'$, where
the representation $V'$ of $\cG_\nabla$ also has $d(V')=0$.)

For example, the adjoint representation ${\mathfrak e}_6$ of $\cG=E_6$
decomposes as a sum of two irreducible representations ${\mathfrak
e}_6 = {\mathfrak f}_4 \oplus V'$ for the differential Galois group
$\cG_\nabla = F_4$, with $\dim V' = 26$. The restriction of $V'$ to
the principal $SL_2$ is the sum $\on{Sym}^{16} \oplus \on{Sym}^8$. The
restriction of $V'$ to $\ch$ is the sum $\chi_1 \oplus \chi_2 \oplus
W_1 \oplus W_2$, with the $\chi_i$ non-trivial of order $3$. Hence
$\on{Irr}_\infty(V') = \dim V'{}^{I_0} = 2, \dim V'{}^{I_\infty} = 0$,
and $d(V') = 0$. 

Similarly, for the non-trivial irreducible representations $V =
\on{Sym}^n$ of $\cG = SL_2$, we find that $d(V) = (n -1)/2$ when $n$
is odd, $d(V) = (n - 2)/2$ when $n$ is congruent to $2$ (mod $4$), and
$d(V) = (n - 4)/2$ when $n$ is divisible by $4$, in agreement with the
results of Section \ref{sl2}. In particular, $d(V) > 0$ whenever $n >
4$. For the spin representation $V$ of dimension $2^n$ of $\cG =
Spin_{2n+1}$, we find that $d(V) > 0$ for all $n > 7$.

\bigskip

We now investigate the difference
\begin{align*}
d(V) &= \on{Irr}_\infty - \dim V^{I_0} - \dim V^{I_\infty} \\
&= \sum_j m(W_j) - \sum_{k \geq 0} m(k) - m(\chi_0)
\end{align*}
further, assuming that $V$ is irreducible and non-trivial. This
argument, which was shown to us by Mark Reeder, breaks into two cases,
depending on whether $\ep=+1$ or $-1$ on $V$.

If $\ep=+1$ on $V$, we may assume (by passing to a quotient that acts
faithfully on $V$) that $\ep=1$ in $\cG$. Then $n^h = 1$ and $\ch$ is
a semi-direct product $\langle n \rangle \ltimes \cS$. In this case,
the map $SL_2 \to \cG$ factors through the quotient $PGL_2$.

We now count the dimension of the span of invariants for $\langle n
\rangle$ on $V$. Since $\langle n \rangle$ is a subgroup of $\ch$
which fixes a line in each $W_j$, we find that
\begin{align*}
\dim V^{\langle n \rangle} &= \sum m(W_j) + m(\chi_0) \\
&= \on{Irr}_\infty + m(\chi_0).
\end{align*}
Since $n$ is conjugate to the element $n' = \rho(e^{2\pi i/h})$, which
lies in the maximal torus $A$ of the principal $PGL_2$, we have
$$
\dim V^{\langle n \rangle} = \sum_{k \geq 0} m(k) \cdot \# \{
\text{weights of $A$ on $\on{Sym}^k$ with $a \equiv 0$ (mod $h$)} \}.
$$
if $m(k) \neq 0$, then $k$ is even and the weight $a=0$ occurs once in
the irreducible representation $\on{Sym}^k$. Since the weights $a$ and
$-a$ occur with the same multiplicity in $V$, we find that
\begin{align*}
\dim V^{\langle n \rangle} &= \sum_{k \geq 0} m(k) + 2 \# \{
\text{weights $a>0$ of $A$ on $V$ with $a \equiv 0$ (mod $h$)} \} \\
&= \dim V^{I_0} + 2 \#  \{
\text{weights $a>0$ of $A$ on $V$ with $a \equiv 0$ (mod $h$)} \}.
\end{align*}

Hence, when $\ep=+1$ we find that
\begin{align*}
d(V) &= \on{Irr}_\infty - \dim V^{I_0} - m(\chi_0) \\
&= 2( \# \{ \text{weights $a>0$ of $A$ on $V$ with $a \equiv 0$ (mod
  $h$)} \} - m(\chi_0)).
\end{align*}
The fact that $d(V)$ is even is consistent with the fact that when $V$
is self-dual and $\ep|_V = +1$, then $V$ is orthogonal. Hence the
first cohomology group $H^1(\pone,j_{!*} \F(V))$ is symplectic.

If the highest weight $\check\la$ of $V$ satisfies
$$
\langle \check\la,\rho \rangle \leq h-1,
$$
then there are no weights $a>0$ with $a \equiv 0$ (mod $h$). Since
$d(V) \geq 0$, this forces $m(\chi_0) = 0 $ and $d(V)=0$. This is what
happens for the adjoint representation.

A similar analysis when $\ep|_V=-1$ gives the formula
\begin{align*}
d(V) &= \# \{ \text{weights $2k+1$ for the torus of $SL_2$ on $V$
  with $k \equiv 0$ (mod $h$) and $k \neq 0$} \} \\ &- m(\chi_1),
\end{align*}
where $\chi_1$ is the character of $\ch$ which is trivial on $\cS$ and
maps $n$ to $e^{\pi i/h}$. if $\ep=n^h$ lies in $\cS$, then
$m(\chi-1)=0$, as all weights of $\cS$ on $V$ are non-trivial.

\section{Nearby connections}    \label{nearby}

There are several connections closely related to the rigid irregular
connection $\nabla$ that we have studied in this paper. First, there
is the connection
$$
\nabla_1 = d + N \frac{dt}{t} 
$$
which has regular singularities at both $t=0$ and $t = \infty$. The
monodromy of $\nabla_1$ is generated by the principal unipotent
element $u = \exp(2 \pi i N)$ and the differential Galois group is the
additive group $\exp(zN)$ of dimension $1$ in $\cG$.

A second related connection is the canonical extension of our local
differential equation at infinity, defined by Katz [K1, 2.4]. This is
the connection
$$
\nabla_2 = d + \left( N - \frac{1}{h}\rho \right) \frac{dt}{t} + E dt.
$$
If we pass to the ramified extension given by $u^h = s = t^{-1}$, and
make a gauge transformation by $g = \rho(u)$ (assuming that $\rho$ is
a co-character), the connection $\nabla_2$ becomes equivalent to the
connection
$$
d - h (N+E) \frac{du}{u^2}.
$$
This is regular at the unique point lying above $t = 0$. The
differential Galois group of $\nabla_2$ and its inertia group at
infinity are both isomorphic to $\ch$. The inertia group at zero is
cyclic, of order $h$ or $2h$.

Finally, we have the following generalization of $\nabla$ for the
exceptional groups $\cG$ of types $G_2$, $F_4$, $E_6$, $E_7$, and
$E_8$, which was suggested by the treatment of nilpotent elements and
regular classes in the Weyl group in \cite{Sp}, Section 9. We thank
Mark Reeder for bringing this argument to our attention.  Let
$$
\varphi': \sw_2 \rightarrow \cg
$$
be a subregular $\sw_2$ in the Lie algebra. The Dynkin labels of the
semi-simple element $\varphi'(h)$ are equal to $2$ on all of the
vertices of the diagram for $\cG$, except the vertex corresponding to
the (unique) root with the highest multiplicity $m$ in the highest
root, where the label is $0$.

Let $d = h - m$. Then the highest eigenspace $\cg[2d-2]$ for
$\varphi'(h)$ on $\cg$ has dimension equal to $1$. Let $E'$ be a
basis, and let $N' = \varphi'(f)$ in $\cg[-2]$. Then Springer shows
that the element $N'+E'$ is both regular and semi-simple in
$\cg$. This element is normalized by the semi-simple element $n' =
\varphi'(h(e^{\pi i/d}))$, which has order either $d$ or $2d$ in
$\cG$.

Let $M$ be the maximal torus which centralizes $N' + E'$. Then the image
of $n'$ in the Weyl group of $M$ is a regular element of order $d$, in
the sense of Springer \cite{Sp}. It has $r+2$ free orbits on the roots
of $\cg$, where $r$ is the rank of $\cG$. We tabulate these numerical
invariants for our groups below, as well as the characteristic
polynomial $F$ of $n'$ acting on the character group of $M$, as a
product of cyclotomic polynomials $F(n)$.

\begin{center}
\begin{tabular}{l|l|l|l|l}
$\cG$  &    $m$   &   $d$   &   $r+2$  &  $F$ \\
\hline
$G_2$  &  3   &     3    &     4    &     $F(3)$ \\
\hline
$F_4$  &  4   &     8    &     6    &     $F(8)$ \\
\hline
$E_6$  &  3   &     9    &     8    &     $F(9)$ \\
\hline
$E_7$  &  4   &    14    &     9    &     $F(14)\cdot F(2)$ \\
\hline
$E_8$  &  6   &    24    &    10    &     $F(24)$
\end{tabular}
\end{center}

\bigskip

Define the subregular analog of $\nabla$ as follows:
$$
\nabla' = d + N' \frac{dt}{t} + E' dt.
$$
Then $\nabla'$ has a regular singularity at $t=0$, with monodromy the
subregular unipotent element $u' = \exp(-2 \pi i N')$. It has an
irregular singularity at $\infty$, with slope $1/d$ and local inertia
group $\langle\cS',n'\rangle$, where $\cS'$ is the subtorus of $M$ on
which $n'$ acts by the primitive $d^{\on{th}}$ roots of unity.

The connection $\nabla'$ is rigid, with differential Galois group
given by the following table

\bigskip

\begin{center}
\begin{tabular}{l|l}
$\cG$   &    $G_{\nabla'}$ \\
\hline
$G_2$   &    $SL_3$ \\
\hline
$F_4$   &   $Spin_9$ \\
\hline
$E_6$   &   $E_6$ \\
\hline
$E_7$   &   $E_7$ \\
\hline
$E_8$   &   $E_8$
\end{tabular}
\end{center}

\bigskip

In the first two cases, we note that a subregular $SL_2$ in $\cG$ is a
regular $SL_2$ in $\cG_{\nabla'}$. Hence the connection $\nabla'$ on
$\cG$ is obtained from the rigid connection $\nabla$ for the
differential Galois group.  In the three other cases, the connection
$\nabla'$ is new; the third gives us a connection with differential
Galois group of type $E_6$. In particular, the differential Galois
group of ${\mathbb G}_m$ has any simple exceptional group as a rigid
quotient.

We get two further rigid connections with differential
Galois group $E_8$ from the two nilpotent classes that
Springer lists in \cite{Sp}, Table 11. These have slopes $1/20$ and
$1/15$ at $\infty$ respectively. Since there is a misprint in that
table, we give the Dynkin labeling of these nilpotent classes below.

$$
\xymatrix{
2 \ar@{-}[r] & 2 \ar@{-}[r] & 0 \ar@{-}[r] & 2 \ar@{-}[r] & 0
\ar@{-}[r] & 2 \ar@{-}[r] & 2 \\
& & & & 2 \ar@{-}[u]}
$$

$$
\xymatrix{
2 \ar@{-}[r] & 0 \ar@{-}[r] & 2 \ar@{-}[r] & 0 \ar@{-}[r] & 2
\ar@{-}[r] & 0 \ar@{-}[r] & 2 \\
& & & & 0 \ar@{-}[u]}
$$
                            
In both cases, the action of the corresponding regular element in the
Weyl group on the character group of the torus is by the primitive
$20^{\on{th}}$ and $15^{\on{th}}$ roots of unity, respectively.

\section{Example of the geometric Langlands correspondence with wild
  ramification}    \label{comments}

The Langlands correspondence for function fields has a counterpart for
complex algebraic curves, called the geometric Langlands
correspondence. As in the classical setting, there is a local and a
global pictures. For simplicity we will restrict ourselves to the case
when $G$ is a simple connected simply-connected algebraic group, so
that $\cG$ is a group of adjoint type.

The local geometric Langlands correspondence has been proposed in
\cite{FG} (see also \cite{Fr} for an exposition). According to
this proposal, to each ``local Langlands parameter'' $\sigma$, which
is a $\cG$-bundle with a connection on $D^\times = \on{Spec}
\C\ppart$, there should correspond a category $\CC_\sigma$ equipped
with an action of the formal loop group $G\ppart$. This correspondence
should be viewed as a ``categorification'' of the local Langlands
correspondence for the group $G(F)$, where $F$ is a local
non-archimedian field, $F = {\mathbb F}_q\ppart$. This means that we
expect that the Grothendieck group of the category $\CC_\sigma$,
equipped with an action of $G\ppart$, attached to a Langlands
parameter $\sigma$, should ``look like'' an
irreducible smooth representation $\pi$ of $G(F)$ attached to an
$\ell$-adic representation of the Weil group of $F$ with the same
properties as $\sigma$.

In particular, an object of $\CC_\sigma$ should correspond to a vector
in the representation $\pi$. Thus, the analogue of the subspace
$\pi^{(K,\Psi)} \subset \pi$ of $(K,\Psi)$-invariant vectors in $\pi$
(that is, vectors that transform via a character $\Psi$ under the
action of a subgroup $K$ of $G(F)$) should be the category
$\CC_\sigma^{(K,\Psi)}$ of $(K,\Psi)$-equivariant objects of
$\CC_\sigma$.

Denote by $\Loc_{\cG}(D^\times)$ the set of isomorphism classes of
flat $\cG$-bundles on $D^\times$. In \cite{FG}, a candidate for the
category $\CC_\sigma$ has been proposed for any $\sigma \in
\Loc_{\cG}(D^\times)$. Namely, consider the category
$\ghat_{\on{crit}}\on{-mod}$ of discrete modules over the affine
Kac--Moody algebra $\ghat$ of critical level. The center of this
category is isomorphic to the algebra of functions on the space
$\on{Op}_{\cG}(D^\times)$ of $\cG$-opers on $D^\times$ (see
\cite{Fr}). Hence for each point $\chi \in \on{Op}_{\cG}(D^\times)$ we
have the category $\ghat_{\on{crit}}\on{-mod}_\chi$ of those modules
on which the center acts according to the character corresponding to
$\chi$.

Consider the forgetful map $p: \on{Op}_{\cG}(D^\times) \to
\Loc_{\cG}(D^\times)$. It was proved in \cite{FZ} that this map is
surjective. Given a local Langlands parameter $\sigma \in
\Loc_{\cG}(D^\times)$, choose $\chi \in p^{-1}(\sigma)$. According to
the proposal of \cite{FG}, the sought-after category $\CC_\sigma$
should be equivalent to $\ghat_{\on{crit}}\on{-mod}_\chi$ (these
categories should therefore be equivalent to each other for all $\chi
\in p^{-1}(\sigma)$). In particular, $\CC_\sigma^{(K,\Psi)}$ should be
equivalent to the category $\ghat_{\on{crit}}\on{-mod}_\chi^{(K,\Psi)}$
of $(K,\Psi)$-equivariant objects in
$\ghat_{\on{crit}}\on{-mod}_\chi$.

Next, we consider the global geometric Langlands correspondence (see
\cite{FG,F:rev,Fr} for more details). Let $X$ be a smooth projective
curve over $\C$. The Langlands parameters $\sigma$ are now
$\cG$-bundles on $X$ with connections which are allowed to have poles
at some points $x_1,\ldots,x_N$. Let $\Bun_{G,(x_i)}$ be the moduli
stack of $G$-bundles on $X$ with the full level structure at each
point $x_i, i=1,\ldots,N$ (that is, a trivialization of the $G$-bundle
on the formal disc at $x_i$). The global correspondence should assign
to $\sigma$ the category $\on{Aut}_{\sigma,(x_i)}$ of Hecke
eigensheaves on $\Bun_{G,(x_i)}$ with eigenvalue $\sigma$.

The compatibility between the local and global correspondence should
be that \cite{FG}
\begin{equation}    \label{eq}
\on{Aut}_{\sigma,(x_i)} \simeq \bigotimes_{i=1}^N \;
\CC_{\sigma_{x_i}},
\end{equation}
where $\sigma_{x_i}$ is the restriction of $\sigma$ to the punctured
disc at $x_i$. This equivalence should give rise to an equivalence of
the equivariant categories. Let $K_{x_i}$ be a subgroup of
$G(\OO_{x_i})$, where $\OO_{x_i}$ is the completed local ring at
$x_i$, and $\Psi_{x_i}$ be its character. Then the equivariant
category of $\on{Aut}_{\sigma,(x_i)}$ is the category
$\on{Aut}^{(K_{x_i},\Psi_{x_i})}_{\sigma,(x_i)}$ of Hecke eigensheaves
on $\Bun_{G,(x_i)}$ with eigenvalue $\sigma$ which are
$(K_{x_i},\Psi_{x_i})$-equivariant. We should have an equivalence of
categories
\begin{equation}    \label{eq eq}
\on{Aut}^{(K_{x_i},\Psi_{x_i})}_{\sigma,(x_i)} \simeq 
\bigotimes_{i=1}^N \; \CC_{\sigma_{x_i}}^{(K_{x_i},\Psi_{x_i})}.
\end{equation}

Now suppose that $\sigma$ comes from an oper $\chi$ which is regular
on $X \bs \{ x_1,\ldots,x_N \}$. Let $\chi_{x_i}$ be the restriction
of $\chi$ to $D_{x_i}^\times$. Then we have the categories
$\ghat_{\on{crit}}\on{-mod}_{\chi_{x_i}}$ and
$\ghat_{\on{crit}}\on{-mod}_{\chi_{x_i}}^{(K_{x_i},\Psi_{x_i})}$,
which we expect to be equivalent to $\CC_{\sigma_{x_i}}$ and
$\CC_{\sigma_{x_i}}^{(K_{x_i},\Psi_{x_i})}$, respectively. The stack
$\Bun_{G,(x_i)}$ may be represented as a double quotient
$$
G(F) \bs \prod_{i=1}^N G(F_{x_i})/\prod_{i=1}^N G(\OO_{x_i}),
$$
where $F = \C(X)$ is the field of rational functions on $X$ and
$F_{x_i}$ is its completion at $x_i$. A loop group version of the
localization functor of Beilinson--Bernstein gives rise to the
functors
\begin{equation}    \label{Delta}
\Delta: \bigotimes_{i=1}^N \; \ghat_{\on{crit}}\on{-mod}_{\chi_{x_i}}
\to \on{Aut}_{\sigma,(x_i)},
\end{equation}
\begin{equation}    \label{Delta1}
\Delta^{(K_{x_i},\Psi_{x_i})}: \bigotimes_{i=1}^N \;
\ghat_{\on{crit}}\on{-mod}_{\chi_{x_i}}^{(K_{x_i},\Psi_{x_i})} \to
\on{Aut}^{(K_{x_i},\Psi_{x_i})}_{\sigma,(x_i)},
\end{equation}
and it is expected \cite{FG} that these functors give rise to the
equivalences \eqref{eq} and \eqref{eq eq}, respectively. This is a
generalization of the construction of Beilinson and Drinfeld in the
unramified case \cite{BD-H}.

\medskip

We now apply this to our situation, which may be viewed as the
simplest example of the geometric Langlands correspondence with wild
ramification (i.e., connections admitting an irregular singularity). We
note that wild ramification has been studied by E. Witten in the
context of $S$-duality of supersymmetric Yang--Mills theory \cite{Wi}.

\medskip

Let $\chi$ be our $\cG$-oper on $\pone$ with poles at the points $0$
and $\infty$. Let $\sigma$ denote the corresponding flat
$\cG$-bundle. We choose $K_0$ to be the Iwahori subgroup $I$, with
$\Psi_0$ the trivial character (we will therefore omit it in the
formulas below), and $K_\infty$ to be its radical $I^0$. Note that
$I^0/[I^0,I^0] \simeq ({\mathbb G}_a)^{\on{rank}(G)+1}$. We choose a
non-degenerate additive character $\Psi$ of $I^0$ as our
$\Psi_\infty$. Thus, we have the global categories
$\on{Aut}_{\sigma,(0,\infty)}$ and
$\on{Aut}^{I,(I^0,\Psi)}_{\sigma,(0,\infty)}$ on
$\Bun_{G,(0,\infty)}$.

According to Section \ref{automorphic}, there is a unique automorphic
representation of $G({\mathbb A})$ corresponding to an $\ell$-adic
analogue of our oper. Moreover, the only ramified local factors in
this representation are situated at $0$ and $\infty$. The former is
the Steinberg representation, whose space of Iwahori invariant vectors
is one-dimensional. The latter is the simple supercuspidal
representation constructed in \cite{GR}. Its space of
$(I^0,\Psi)$-invariant vectors is also one-dimensional. We recall that
the geometric analogue of the space of invariant vectors is the
corresponding equivariant category. Hence the geometric counterpart of
the multiplicity one statement of Section \ref{automorphic} is the
statement (conjecture) that {\em the category
$\on{Aut}^{I,(I^0,\Psi)}_{\sigma,(0,\infty)}$ has a unique
non-zero irreducible object}. (Here and below ``unique'' means
``unique up to an isomorphism.'')

\medskip

The compatibility of the local and global correspondences gives us a
way to construct this object. Namely, we have two local categories
$\ghat_{\on{crit}}\on{-mod}_{\chi_0}^I$ and
$\ghat_{\on{crit}}\on{-mod}_{\chi_\infty}^{(I^0,\Psi)}$ attached to
the points $0$ and $\infty$, respectively. The oper $\chi_0$ on
$D_0^\times$ has regular singularity and regular unipotent
monodromy. Using the results of \cite{FG,Fr}, one can show that the
category $\ghat_{\on{crit}}\on{-mod}_{\chi_0}^I$ has a unique non-zero
irreducible object ${\mathbb M}_{-\rho}(\chi_0)$, which is constructed
as follows. It is the quotient of the Verma module
$$
{\mathbb M}_{-\rho} = \on{Ind}_{\on{Lie}(I) \oplus \C {\mathbf
1}}^{\ghat_{\on{crit}}}(\C_{-\rho})
$$
over $\ghat_{\on{crit}}$ with highest weight $-\rho$, by the image of
the maximal ideal in the center corresponding to the central character
$\chi_0$.

On the other hand, $\chi_\infty$ has irregular singularity with the
slope $1/h$. To construct an object of the category of
$\ghat_{\on{crit}}\on{-mod}_{\chi_\infty}^{(I^0,\Psi)}$, we imitate
the construction of \cite{GR} (see Section \ref{automorphic}). Define
the affine Whittaker module
$$
{\mathbb W}_{\Psi} = \on{Ind}_{\on{Lie}(I^0) \oplus \C {\mathbf
1}}^{\ghat_{\on{crit}}}(\Psi)
$$
over $\ghat_{\on{crit}}$ (here we denote by the same symbol $\Psi$ the
character of the Lie algebra $\on{Lie}(I^0)$ corresponding to the
above character $\Psi$ of the group $I^0$). Let ${\mathbb
W}_\Psi(\chi_\infty)$ be the quotient of ${\mathbb W}_\Psi$ by the
image of the maximal ideal in the center corresponding to the central
character $\chi_\infty$. By construction, it is an
$(I^0,\Psi)$-equivariant $\ghat_{\on{crit}}$-module and hence it is
indeed an object of our local category
$\ghat_{\on{crit}}\on{-mod}_{\chi_\infty}^{(I^0,\Psi)}$. Applying the
localization functor $\Delta^{I,(I,\Psi)}$ of \eqref{Delta1} to
${\mathbb M}_{-\rho}(\chi_0) \otimes {\mathbb W}_\Psi(\chi_\infty)$,
we obtain an object of the category
$\on{Aut}^{I,(I^0,\Psi)}_{\sigma,(0,\infty)}$. It is
natural to conjecture that this is the unique non-zero irreducible
object of this category.

One can show (see \cite{FF:kdv}, Lemma 5) that the image of the center
of the completed enveloping algebra in $\on{End} {\mathbb W}_\Psi$ is
the algebra of functions on the space of opers which have
representatives of the form
$$
d - N \frac{ds}{s} - E \frac{ds}{s^2} + {\mathbf v} \frac{ds}{s},
$$
where ${\mathbf v} \in \bb[[s]]$. Since our oper $\chi_\infty$ belongs
to this space, we obtain that the quotient ${\mathbb
W}_\Psi(\chi_\infty)$ is non-zero. This provides supporting evidence
for the above conjecture describing an example of
the geometric Langlands correspondence with wild ramification.

\end{document}